\documentclass[preprint,10pt]{elsarticle}

\usepackage{amssymb}
\usepackage{amsmath}
\usepackage{lineno}
\usepackage{algorithmic}
\usepackage{algorithm2e}
\usepackage{subcaption}
\usepackage{relsize}
\usepackage{graphicx}
\usepackage[utf8x]{inputenc}
\usepackage[title]{appendix}
\journal{arXiv.org}

\newcommand{\dext}{\mathrm{d}}
\renewcommand{\vec}[1]{\mathbf{#1}}
\newcommand{\norm}[1]{\left\lVert#1\right\rVert}

\newcommand{\cpsi}{c_{\psi}}
\newcommand{\spsi}{s_{\psi}}
\newcommand{\ctheta}{c_{\theta}}
\newcommand{\stheta}{s_{\theta}}
\newcommand{\cgamma}{c_{\gamma}}
\newcommand{\sgamma}{s_{\gamma}}
\newcommand{\cphi}{c_{\phi}}
\newcommand{\sphi}{s_{\phi}}

\begin{document}
\begin{frontmatter}

\title{Integrable cross-field generation based on imposed singularity configuration \\
 -- the 2D manifold case --}
\author{Jovana Jezdimirovi\'c, Alexandre Chemin, Jean-Fran\c{c}ois Remacle}
\address{Universit\' e catholique de Louvain, Louvain la Neuve, Belgium 
jean-francois.remacle@uclouvain.be}

\begin{abstract}
This work presents the mathematical foundations for the generation of integrable cross-field on 2D manifolds based on user-imposed singularity configuration.
In this paper, we either use singularities that appear naturally, \emph{e.g.}, 
by solving a non-linear problem, or use as an input user-defined singularity pattern, 
possibly with high valence singularities that typically do not appear in cross-field computations. This singularity set is under the constraint of Abel-Jacobi's conditions for valid singularity configurations.
The main contribution of the paper is the development of a formulation that allows computing an integrable isotropic 2D cross-field from a given set of singularities through the resolution of only two linear PDEs. 
To address the issue of possible suboptimal singularities' distribution, we also present the mathematical setting for the generation of an integrable anisotropic 2D cross-field based on a user-imposed singularity pattern.
The developed formulations support both an isotropic and an anisotropic block-structured quad mesh generation.
\end{abstract}


\end{frontmatter}

\section{Introduction and related work}
\label{S:Intro}

Numerous methods for surface parametrization/representation have been developed for a large number of applications  \cite{bommes2012, campen2017, Floater:2005}. In cases when a shape exhibits complex topological or geometrical characteristics, it is necessary to split it into simple partitions to obtain a quad mesh. The special case of partitioning it into a simply connected network of conformal quadrilateral partitions is called the \textit{quad layout} \cite{campen2014}. The latter manner of surface representation is a subject of great interest in meshing and computer graphics communities, due to providing a wide range of benefits \cite{bommes2012, campen2017, Shepherd:2020}. Nevertheless, these advantages come with the high price of dealing with complex and time-consuming algorithms \cite{Shepherd:2020}. 

Among the developed methods for the quad layout generation, a general distinction can be made among the ones which are: computing a seamless global parametrization of the domain where integer iso-values of the parameter fields form the sides \cite{campen2015, bommes2013, ray2006}, using Riemann geometry \cite{Chen:2019, Lei:2020, zheng2021quadrilateral}, or like in our case, constructing a \textit{cross-field} structure that will guide the integral lines emanating from \textit{singularities} \cite{kowalski2013pde, Fogg:2015,  jezdimirovic2019,  Jezdimirovic:2021, ray2008,  viertel2019approach}. 

Although leaning on heterogeneous approaches, all the above-mentioned methods share the common challenge: dealing with the inevitable singularity configuration. A \textit{singularity} appears where a cross-field vanishes and it represents an irregular vertex of a quad layout/quad mesh \cite{beaufort2017}, \emph{i.e.}, a vertex which doesn't have exactly four adjacent quadrilaterals.
The singular configuration is constrained by the Euler characteristic $\chi$, which is a topological invariant of a surface. Moreover, a suboptimal number or location of singularities can have severe consequences: causing undesirable  thin partitions, large distortion, not an adequate number and/or tangential crossings of separatrices as well as limit cycles \cite{Shepherd:2020, bommes2013, viertel2019approach}.  

\begin{figure*} 
\begin{center}
  \includegraphics[width=0.99\linewidth]{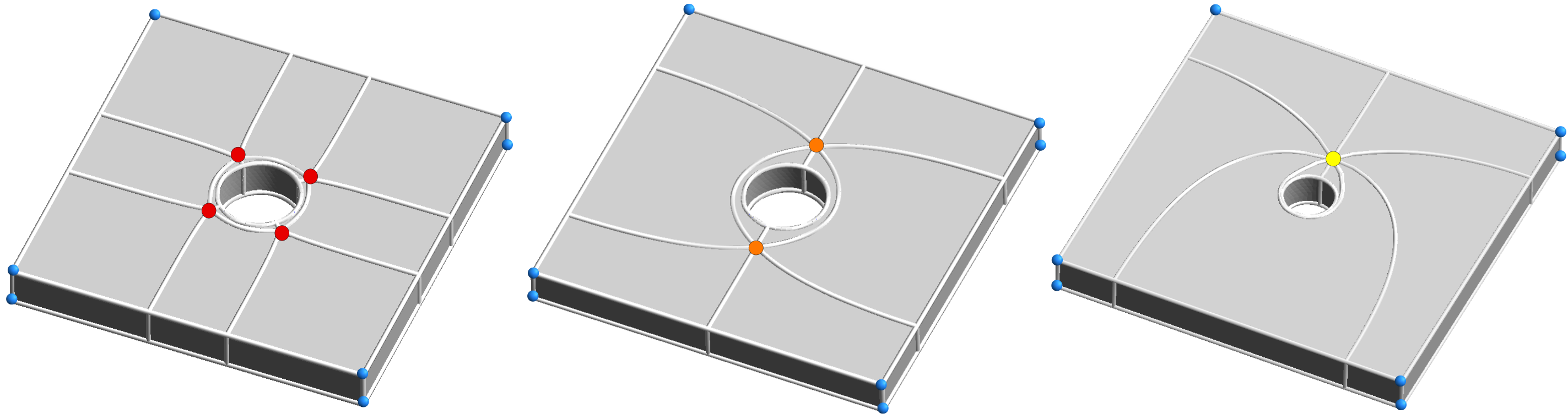}
  \caption{Three quad layouts of a simple domain. Singularities of valence $3$ are colored in blue, valence $5$ in red, valence $6$ in orange, and valence $8$ in yellow.}
  \label{fig:teaser}
\end{center}
\end{figure*}

Cross-field guided methods can be very useful and flexible but they typically lack direct control over the positions of the singularities and the structures of the quad layout \cite{gu2020computational}. Our cross-field formulation, with mathematical foundations detailed in Section~\ref{S:Cross-field}, offers a contribution to this issue through the concept of user-imposed singularity configuration in order to gain direct control over their number, location, and valence (number of adjacent quadrilaterals).
The user is entitled to use either naturally appearing singularities, obtained by solving a non-linear problem \cite{jezdimirovic2019, viertel2019approach, vaxman2016directional, Hertzmann:2000}, using globally optimal direction fields \cite{knoppel2013globally}, or to impose its own singularity configuration, possibly with high valences, as illustrated in Fig.~\ref{fig:teaser}. It is important to note that the choice of singularity pattern is not arbitrary, though. Moreover, it is under the direct constraint of \textit{Abel-Jacobi theory} \cite{Chen:2019, Lei:2020, zheng2021quadrilateral} for valid singularity configurations. 
Here, the singularity configuration is taken as an input and an integrable isotropic cross-field is computed by solving only two linear systems, Section~\ref{S:Isotropic}. Computation of the scalar field $H$ used for this cross-field generation bears some resemblance to the one developed for unstructured mesh generation on planar and curved surface domains in \cite{bunin2008continuum}. 
Finally, the preliminary results of the developed cross-field formulation for an isotropic block-structured quad mesh generation are outlined using the 3-step pipeline \cite{Jezdimirovic:2021} in Section~\ref{S:Preliminary}. 

Computing only one scalar field $H$ (a metric that is flat except at singularities) imposes a strict constraint on  singularities' placement, \emph{i.e.}, fulfilling all Abel-Jacobi conditions. In practice, imposing suboptimal distribution of singularities may lead to not obtaining boundary-aligned cross-field, disabling an isotropic quad mesh generation, Section~\ref{sec:nonQuadMeshableSingConfig} and ~\ref{sec:SubOptimalSing}. To bypass this issue, we develop a new cross-field formulation on the imposed singularity configuration, which considers the integrability, while relaxing the condition on isotropic scaling of crosses' branches. Here, two independent metrics $H_1$ and $H_2$ are computed instead of only one as in the Abel-Jacobi framework, enabling an integrable 2D cross-field generation with anisotropic scaling, Section~\ref{S:Anisotropic}.  

Lastly, final remarks and some of the potential applications are discussed in Section~\ref{S:Conclusion}.

\section{Cross-field computation on prescribed singularity configuration}
\label{S:Cross-field}

We define a 2D cross $\mathbf{c}$ as a set of $2$ unit coplanar orthogonal vectors and their opposite, \emph{i.e.}, $$\mathbf{c}=\{ \mathbf{u} ,\mathbf{v} ,-\mathbf{u} ,-\mathbf{v}\}$$ with $\{\mathbf u . \mathbf v=0$, $|\mathbf{u}|=|\mathbf{v}|=1\}$ and $\mathbf u, \mathbf v$ are coplanar. These vectors are called cross' branches.

A 2D cross-field $\mathcal{C}_\mathcal{M}$ on a 2D manifold $\mathcal M$, now, is a map $\mathcal{C}_\mathcal{M}: \mathbf{X} \in \mathcal{M} \to \mathbf{c}(\mathbf{X})$,
and the standard approach to compute a smooth boundary-aligned cross-field is to minimize the Dirichlet energy:
\begin{equation}
\displaystyle \min_{\mathcal{C}_\mathcal{M}} \int_\mathcal{M} \norm{\nabla \mathcal{C}_\mathcal{M}}^2 
\label{eq:dirichlet_energy}
\end{equation}
subject to the boundary condition $\mathbf{c(\mathbf{X})} = \mathbf{g}(\mathbf{X})$ on $\partial \mathcal{M}$,
where $\mathbf{g}$ is a given function.

The classical boundary condition for cross-field computation is that $\forall \mathbf P \in \partial \mathcal M$, with $\mathbf T (\mathbf P)$ a unit tangent vector to $\mathcal M$ at $\mathbf P$, one branch of $\mathbf c (\mathbf P)$ has to be colinear to $\mathbf T (\mathbf P)$. In the general case, there exists no smooth cross-field matching this boundary condition. The cross-field will present a finite number of singularities $\mathbf S_j$, located at $\mathbf X_j$ and of index $k_j$, related to the concept of valence as $k_j = 4-\mathrm{valence}(\mathbf S_j).$ \\
We define a singularity configuration as the set $$\mathcal S = \{ \mathbf S_j, j\in[|1,N |], N\in\mathbb{Z}\}.$$

In the upcoming section, a method to compute a cross-field $\mathcal C_\mathcal{M}$ matching a given singularity configuration $\mathcal S$ is developed. In other words, we are looking for $\mathcal C_\mathcal{M}$ such as:
\begin{equation}
\left\{
\begin{array}{l}
\text{$\boldsymbol{\cdot}$ if $\vec X$ belongs to $\partial \mathcal M$, at least one branch of} 
\text{ $\mathcal C_\mathcal{M}(\vec X)$ is tangent to $\partial \mathcal M$,}  \\
\text{$\boldsymbol{\cdot}$ singularities of $\mathcal C_\mathcal{M}$ are matching the given $\mathcal S$}\\
\text{ (the same number, location, and indices).} 
\end{array}
\right.
\label{eq:cond}
\end{equation}

Before developing the method to compute such a cross-field, a few operators on the 2D manifold have to be defined.

\subsection{Curvature and Levi-Civita connection on the 2D manifold}
\label{sec:def}
  
Let $E^3$ be the Euclidean space equipped with a Cartesian coordinates system $\{x^i, i=1,2,3\}$, and 
$\mathcal{M}$ be an oriented two-dimensional manifold embedded in $E^3$. We note $\mathbf n(\mathbf X)$ the unit normal to $\mathcal{M}$ at $\mathbf X \in \mathcal{M}$. It is assumed that the normal field $\mathbf n$ is smooth and that the Gaussian curvature $K$ is defined and smoothed on $\mathcal M$.

If $\gamma(s)$ is a curve on $\mathcal{M}$ {\em parametrized by arc length}, the Darboux frame is the orthonormal frame defined by
\begin{align}
\mathbf T(s) &= \gamma'(s)\\
\mathbf n(s) &= \mathbf n(\gamma(s))\\
\mathbf t(s) &= \mathbf n(s) \times \mathbf T(s).
\end{align}
One then has the differential structure
\begin{equation}
\dext \left( \begin{array}{c} \mathbf T\\ \mathbf t\\ \mathbf n\end{array} \right) 
= \left( \begin{array}{ccc} 
0 & \kappa_g & \kappa_n\\ 
-\kappa_g & 0  & \tau_r \\ 
-\kappa_n & -\tau_r & 0 \end{array} \right)
\left( \begin{array}{c} \mathbf T\\ \mathbf t\\ \mathbf n\end{array} \right) \dext s
\label{darboux}
\end{equation}
where $\kappa_g$ is the geodesic curvature of the curve, $\kappa_n$ the normal curvature of the curve, and $\tau_r$ the relative torsion of the curve. $\mathbf T$ is the unit tangent, $\mathbf t$ the tangent normal and $\mathbf n$ the unit normal.

Arbitrary vector fields $V$ and $W\in E^3$ 
can be expressed as 
$$
\mathbf V = V^i \mathbf E_i,
\quad  \quad 
\mathbf W = W^i \mathbf E_i 
$$
in the natural basis vectors $\{\mathbf E_i, i=1,2,3\}$ of this coordinate system, 
and we shall note
$$
< \mathbf V, \mathbf W> = V^i W^j \delta_{ij},
\quad  \quad 
|| \mathbf V || = \sqrt { < \mathbf V, \mathbf V > } 
$$ 
the Euclidean metric and the associated norm for vectors.
The Levi-Civita connection on $E^3$ in Cartesian coordinates 
is trivial (all Christoffel symbols vanish), and one has
$$
\nabla^E_{\mathbf V} \mathbf W 
= (\nabla_{\mathbf V}  W^i) \mathbf E_i. 
$$

The Levi-Civita connection on the Riemannian submanifold $\mathcal{M}$,
now, is not a trivial one.
It is the orthogonal projection of $\nabla^E_{\mathbf V}$
in the tangent bundle 
$T\mathcal{M}$, so that one has
\begin{equation}
\nabla_{\mathbf V} \mathbf W =  P_{T\mathcal{M}}[ \nabla^E_{\mathbf V} \mathbf W]
=  (\nabla_{\mathbf V}  W^i) P_{T\mathcal{M}}[ \mathbf E_i  ]
\label{connection}
\end{equation}
where $P_{T\mathcal{M}}\,: E^3 \mapsto T\mathcal{M}$ 
is the orthogonal projection operator on $T\mathcal M$.

An arbitrary orthonormal local basis $(\mathbf{u}_\mathbf{X},\mathbf{v}_\mathbf{X},\mathbf{n})$ for every $\mathbf{X}\in\mathcal{M}$, can be represented through the Euler angles $(\psi,\gamma,\phi)$ which are $\mathcal{C}^1$ on $\mathcal{M}$, and with the shorthands $s_\phi \equiv \sin \phi$ and $c_\phi \equiv \cos \phi$, as:
\begin{equation}
\begin{array}{ccc}
\mathbf{u}_\mathbf{X} & = & \left(\begin{matrix}- \sphi \spsi \cgamma + \cphi \cpsi\\\sphi \cpsi \cgamma + \spsi \cphi\\\sphi \sgamma\end{matrix}\right),\,  \\
\mathbf{v}_\mathbf{X} & = & \left(\begin{matrix}- \sphi \cpsi - \spsi \cphi \cgamma\\- \sphi \spsi + \cphi \cpsi \cgamma\\\sgamma \cphi\end{matrix}\right),\,  \\
\mathbf n & = & \left(\begin{matrix}\spsi \sgamma\\- \sgamma \cpsi\\\cgamma\end{matrix}\right)
\label{eulerangles}
\end{array}
\end{equation}
in the vector basis of $E^3$.
 
\subsection{Conformal mapping}
\label{sec:Conformal}

We are looking for a conformal mapping
\begin{equation}
  \begin{array}{ccrcl}
    \mathcal{F}&:& \mathcal{P} & \rightarrow & \mathcal{M} \subset E^3\\
    & &  \mathbf{P}=(\xi,\eta) & \mapsto & \mathbf{X}=(x^1,x^2,x^3)
  \end{array}
\end{equation}
where $\mathcal P$ is a parametric space.
As finding $\mathcal F$ right away is a difficult problem,
one focuses instead on finding the $3 \times 2$ jacobian matrix of $\mathcal F$
\begin{equation}
 J(\mathbf P) = ( \partial_\xi \mathcal F(\mathbf P), \partial_\eta \mathcal F(\mathbf P) )
\equiv ( \tilde{\mathbf u}(\mathbf P), \tilde{\mathbf v}(\mathbf P) ),
\end{equation}
where $\tilde{\mathbf u}, \tilde{\mathbf v} \in T\mathcal{M}$ are the columns vectors of $J$. 
The mapping $\mathcal F$ being conformal, 
the columns of $J(\mathbf P)$ have the same norm $L(\mathbf P) \equiv || \tilde{\mathbf u}(\mathbf P) || = || \tilde{\mathbf v}(\mathbf P) ||$ 
and are orthogonal to each other, $\tilde{\mathbf u}(\mathbf P) \cdot \tilde{\mathbf v}(\mathbf P) =0$.
We can also write:
$$
J = L  ( \mathbf u, \mathbf v ),
\quad  \quad 
\mathbf n = \mathbf u \land \mathbf v
$$
where 
\begin{equation}
  \begin{array}{rcl}
    \mathbf u & = & \frac{\tilde{\mathbf u}}{||\tilde{\mathbf u}||}\\
    \mathbf v & = & \frac{\tilde{\mathbf v}}{||\tilde{\mathbf v}||}.
  \end{array}
\end{equation}
Recalling that finding a conformal transformation $\mathcal{F}$ is challenging, we will from now on be looking for the jacobian $J$, \textit{i.e.}, the triplet $(\mathbf u, \mathbf v, L)$.

The triplet $(\mathbf u, \mathbf v, \mathbf n)$ forms a set of $3$ orthonormal basis vectors and can be seen as a rotation of $(\mathbf{u}_\mathbf{X}, \mathbf{v}_\mathbf{X}, \mathbf n)$ among the direction $\mathbf n$. Therefore, a 2D cross ${c}(\mathbf X)\text{, }\mathbf X\in \mathcal{M}$ can be defined with the help of a scalar field $\theta$, where $\mathbf{u} = \mathcal{R}_{\theta,\mathbf{n}}(\mathbf{u}_{\mathbf{X}})$ and $\mathbf{v} = \mathcal{R}_{\theta,\mathbf{n}}(\mathbf{v}_{\mathbf{X}})$, and the local manifold basis $(\mathbf{u}_\mathbf{X}, \mathbf{v}_\mathbf{X}, \mathbf n)$ as:
\begin{equation}
\mathbf u = \ctheta \mathbf{u}_\mathbf{X} + \stheta \mathbf{v}_\mathbf{X}\text{, } \quad \quad 
\mathbf v = -\stheta \mathbf{u}_\mathbf{X} + \ctheta \mathbf{v}_\mathbf{X}.
\end{equation}

By using the Euler angles $(\psi,\gamma,\phi)$ and $\theta$, the triplet $(\mathbf u, \mathbf v, \mathbf n)$ can also be expressed as:
\begin{equation}
\begin{array}{ccc}
\mathbf u & = & \left(\begin{matrix}- s_{\theta+\phi} \spsi \cgamma + c_{\theta+\phi} \cpsi\\s_{\theta+\phi} \cpsi \cgamma + \spsi c_{\theta+\phi}\\s_{\theta+\phi} \sgamma\end{matrix}\right),\, \\
\mathbf v & = & \left(\begin{matrix}- s_{\theta+\phi} \cpsi - \spsi c_{\theta+\phi} \cgamma\\- s_{\theta+\phi} \spsi + c_{\theta+\phi} \cpsi \cgamma\\\sgamma c_{\theta+\phi}\end{matrix}\right),\, \\
\mathbf n & = & \left(\begin{matrix}\spsi \sgamma\\- \sgamma \cpsi\\\cgamma\end{matrix}\right).
\label{euleranglesCross}
\end{array}
\end{equation}
It is important to note that $\mathbf u$ and $\mathbf v$ are the two branches of the cross-field $\mathcal C_\mathcal{M}$ we are looking for. The projection operator $P_{T\mathcal{M}}$ introduced in Eq.~(\ref{connection}) then simply amounts to disregarding the component along $\mathbf n$ of vectors.

For a vector field $\mathbf w$ defined on $\mathcal{M}$,
one can write by derivation of Eq.~(\ref{euleranglesCross})
\begin{equation}
\begin{aligned}
\nabla^E_{\mathbf w} \mathbf u &= \mathbf v \,\nabla_{\mathbf w}(\theta+\phi) 
+ s_{\theta+\phi} \mathbf n \,\nabla_{\mathbf w} \gamma \\
& \quad + (c_\gamma \mathbf v - s_\gamma c_{\theta + \phi} \mathbf n) \nabla_{\mathbf w} \Psi \\
\nabla^E_{\mathbf w} \mathbf v &= -\mathbf u \,\nabla_{\mathbf w}(\theta+\phi) 
+ c_{\theta + \phi} \mathbf n \,\nabla_{\mathbf w} \gamma \\
& \quad + (-c_\gamma \mathbf u + s_\gamma s_{\theta+\phi} \mathbf n) \nabla_{\mathbf w} \Psi \\
\nabla^E_{\mathbf w} \mathbf n &=  - (s_{\theta+\phi} \mathbf u + c_{\theta + \phi} \mathbf v) \nabla_{\mathbf w} \gamma \\
& \quad + s_\gamma ( c_{\theta + \phi} \mathbf u - s_{\theta+\phi} \mathbf v)  \nabla_{\mathbf w} \Psi 
\end{aligned}
\end{equation}
and hence, using Eq.~(\ref{connection}), 
the expression of the covariant derivatives on the submanifold $\mathcal{M}$ is:
\begin{equation}
\begin{aligned}
\nabla_{\mathbf w} \mathbf u &= \mathbf v \,\nabla_{\mathbf w}(\theta+\phi) + c_\gamma \mathbf v  \nabla_{\mathbf w} \Psi \\
\nabla_{\mathbf w} \mathbf v &= -\mathbf u \,\nabla_{\mathbf w}(\theta+\phi) - c_\gamma \mathbf u \nabla_{\mathbf w} \Psi.
\end{aligned}
\end{equation}
This allows writing the Lie bracket
\begin{align}
[\mathbf u , \mathbf v] &= \nabla_{\mathbf u} \mathbf v - \nabla_{\mathbf v} \mathbf u \nonumber\\
&=  -( \mathbf u \nabla_{\mathbf u} (\theta+\phi) + \mathbf v \nabla_{\mathbf v} (\theta+\phi) \nonumber \\ 
& \quad - c_\gamma (  \mathbf u \nabla_{\mathbf u} \psi + \mathbf v \nabla_{\mathbf v} \psi ),
\label{bracketuv}
\end{align}
which will be used in the upcoming section.
    
\section{Integrability condition with isotropic scaling}
\label{S:Isotropic}

The mapping $\mathcal F$, now, 
defines a conformal parametrization of $\mathcal{M}$
if the columns of $J$ commute as vector fields,
\emph{i.e.}, if the differential condition
\begin{equation}
0 = [\tilde{\mathbf{u}},\tilde{\mathbf{v}}] 
= \nabla_{\tilde{\mathbf u}}  \tilde{\mathbf v} -\nabla_{\tilde{\mathbf v}} \tilde{\mathbf u} 
= [ L \mathbf u , L \mathbf v] 
\label{integrability}
\end{equation}
is verified.
Developing the latter expression and posing for convenience $L = e^H$, it becomes
$$
0 = \mathbf v \nabla_{\mathbf u} H - \mathbf u \nabla_{\mathbf v} H + [\mathbf u , \mathbf v],
$$
and then
\begin{equation}
  \left\{
  \begin{array}{crcl}
    \nabla_\mathbf{u} H & = & - < \mathbf v,  [\mathbf u , \mathbf v]  > \\
    \nabla_\mathbf{v} H & = & < \mathbf u,  [\mathbf u , \mathbf v]  >
  \end{array}
  \right.
  \label{intManifold}
\end{equation}
which after the substitution of Eq.~(\ref{bracketuv}) gives
\begin{equation}
  \left\{
  \begin{array}{rrl}
    \nabla_{\mathbf u} H & = &  \nabla_{\mathbf v}\theta\, + \nabla_{\mathbf v}\phi\, + \cgamma \nabla_{\mathbf v}\psi  \\
    -\nabla_{\mathbf v} H & = & \nabla_{\mathbf u}\theta\, + \nabla_{\mathbf u}\phi\, + \cgamma \nabla_{\mathbf u}\psi.
  \end{array}
  \right.
  \label{eq:intManifold2}
\end{equation}

In order to obtain the boundary value problem for $H$, the partial differential equation (\emph{PDE}) governing it will be expressed on $\partial \mathcal{M}$ as well as on the interior of $\mathcal{M}$.

\subsection{$H$ PDE on the boundary}
\label{sec:PDEbndH}
As the boundary $\partial \mathcal{M}$ is represented by curves on $\mathcal{M}$, it is possible to parametrize them by arc length and thus associate for each $\mathbf X \in \partial \mathcal{M}$ a Darboux frame $(\mathbf T (\mathbf X), \mathbf t(\mathbf X), \mathbf n(\mathbf X))$. As we are looking for a cross-field $\mathcal C_\mathcal{M}$ fulfilling conditions~(\ref{eq:cond}), the triplet $(\mathbf u, \mathbf v, \mathbf n)$ can be identified as $(\mathbf T (\mathbf X), \mathbf t(\mathbf X), \mathbf n(\mathbf X))$. One then has:
$$
\begin{array}{lcc}
\partial_s T = \kappa_g \mathbf t + \kappa_n \mathbf n
\equiv \nabla_{\mathbf u}\mathbf u \\
\quad \quad = \mathbf v \nabla_{\mathbf u}\phi + s_\phi \mathbf n \nabla_{\mathbf u}\gamma 
+ (c_\gamma \mathbf v - s_\gamma c_\phi \mathbf n) \nabla_{\mathbf u}\psi 
\end{array}
$$
where from follows
\begin{equation}
  \left\{
  \begin{array}{rcl}
    \kappa_g &=& \nabla_{\mathbf u}\phi + c_\gamma \nabla_{\mathbf u}\psi \\
    \kappa_n &=& s_\phi \nabla_{\mathbf u}\gamma -s_\gamma c_\phi \nabla_{\mathbf u}\psi.
  \end{array}
  \right.
  \label{eq:curv}
\end{equation}
Using Eq.~(\ref{eq:intManifold2}) it becomes:
\begin{equation}
  \nabla_{\mathbf t} H = -\kappa_g,
  \label{eq:bndH}
\end{equation}

the result that matches exactly the one found in the planar case \cite{Jezdimirovic:2021}.

\subsection{$H$ PDE in the smooth region on the interior of $M$}
\label{sec:localEquationInSmoothRegion}

To find the \emph{PDE} governing $H$, let's assume the jacobian $J$ is smooth (and therefore $H$) in a vicinity $\mathcal V$ of $\mathbf X \in M$.

We choose $\mathcal{U} \subset \mathcal V$ such as $\mathbf{X}\in\mathcal{U}$, $\partial\mathcal{U}$ such as unit tangent vector $\mathbf T _0$ to $\partial\mathcal{U}_0$ verifies $\mathbf T_0=\mathbf{v}$, $\mathbf T_1$ to $\partial\mathcal{U}_1$ verifies $\mathbf T_1=\mathbf{u}$, $\mathbf T_2$ to $\partial\mathcal{U}_2$ verifies $\mathbf T_2=-\mathbf{v}$, $\mathbf T_3$ to $\partial\mathcal{U}_3$ verifies $\mathbf T_3=-\mathbf{u}$.

\begin{figure}[h!t]
  \centering
  \includegraphics[width=0.35\textwidth]{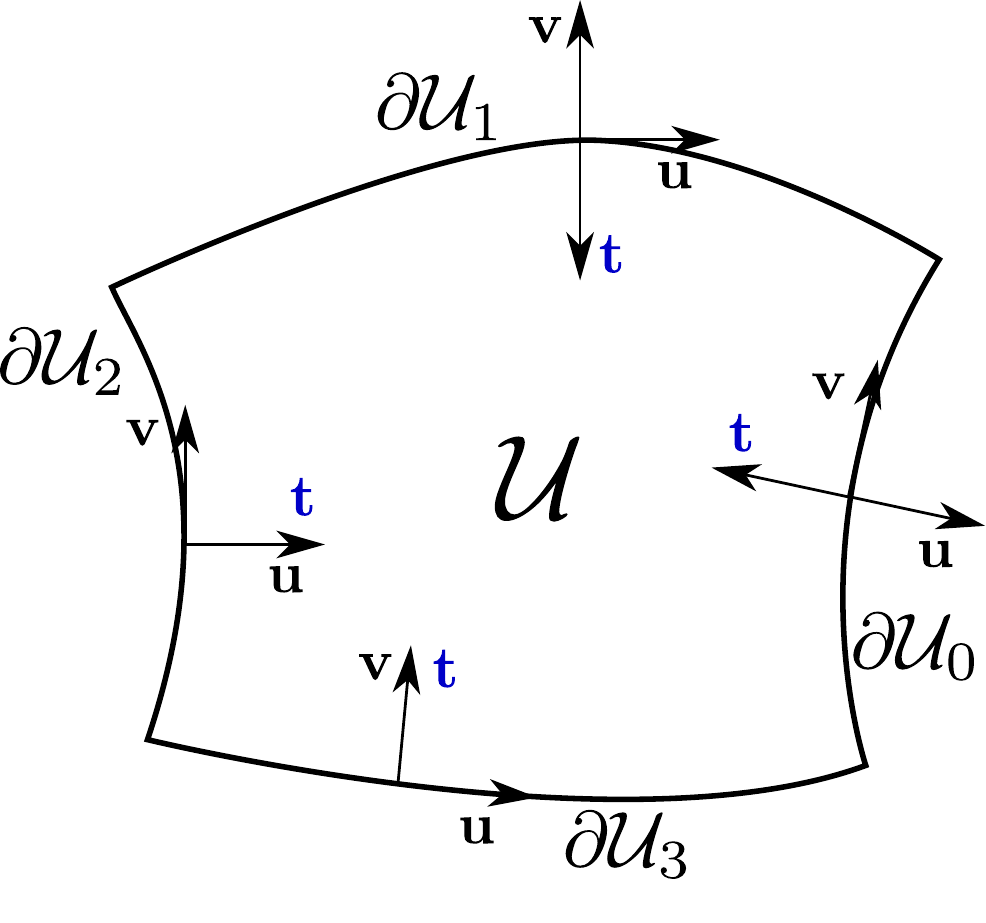}
  \caption{Vicinity of $\mathbf X$ considered.}
  \label{fig:vicinity}
\end{figure}

Thus we have a submanifold $\mathcal U \subset M$ on which $H$ is smooth, and such as $\partial\mathcal{U}=\partial\mathcal U_0\cup\partial\mathcal U_1\cup\partial\mathcal U_2\cup\partial\mathcal U_3$. Darboux frames of $\partial\mathcal U$ (Fig.~\ref{fig:vicinity}) are:
\begin{equation}
  \left\{
  \begin{array}{rcll}
    (\mathbf T, \mathbf t, \mathbf n)&=& (\,\,\,\,\,\mathbf v, -\mathbf u, \mathbf n) &\text{on } \partial\mathcal U_0\\
    (\mathbf T, \mathbf t, \mathbf n)&=& (-\mathbf u, -\mathbf v, \mathbf n) &\text{on } \partial\mathcal U_1\\
    (\mathbf T, \mathbf t, \mathbf n)&=& (-\mathbf v, \,\,\,\,\,\mathbf u, \mathbf n) &\text{on } \partial\mathcal U_2\\
    (\mathbf T, \mathbf t, \mathbf n)&=& (\,\,\,\,\,\mathbf u, \,\,\,\,\,\mathbf v, \mathbf n) &\text{on } \partial\mathcal U_3\\
  \end{array}
  \right.
\end{equation}

For $(\mathbf{\tilde{u}}, \mathbf{\tilde{v}})$ to be a local coordinate system, we recall Eq.~(\ref{eq:bndH}) demonstrated in Section~\ref{sec:PDEbndH}:
\begin{equation}
  \kappa_g = -\nabla_{\mathbf{t}} H\text{, with }\mathbf{t} = \mathbf{n}\land\mathbf{T}
\label{eq:smooth1}
\end{equation}
and the divergence theorem stating that:
\begin{equation}
  \displaystyle\int_{\partial\mathcal{U}}\nabla_\mathbf{t} H = -\displaystyle\int_{\mathcal{U}} \Delta H.
\label{eq:smooth2}
\end{equation}

Applying the Gauss-Bonnet theorem on $\mathcal{U}$ leads to:
$$\int_\mathcal{U} K \,\text{d}\mathcal{U} + \int_{\partial \mathcal{U}} \kappa_g \,\text{d}l + 4\frac{\pi}{2}= 2 \pi \chi(\mathcal{U})$$
where $K$  and $\chi(\mathcal{U})$ are respectively the Gaussian curvature and the Euler characteristic of $\mathcal U$. As $\chi(\mathcal{U})=1$ and using Eq.~(\ref{eq:smooth1}) and (\ref{eq:smooth2}), it becomes:
\begin{equation}
  \int_\mathcal{U}K\,\text{d}\mathcal{U} = -\int_{\mathcal{U}}\Delta H\,\text{d}\mathcal{U}
\end{equation}
which holds for any chosen $\mathcal{U}$. Hence, there is:
\begin{equation}
  \Delta H=-K\text{, if } J \text{ is smooth}.
  \label{eq:smoothH}
\end{equation}

In the general case, it is impossible for $J$ to be smooth everywhere. Indeed, let's assume $\mathcal M$ to be with smooth boundary $\partial \mathcal M$ (\emph{i.e.} with no corners) and of the Euler characteristic $\chi(\mathcal M)=1$. If we assume $J$ is smooth everywhere, it becomes:
\begin{equation}
  \left\{
  \begin{array}{lll}
    \int_\mathcal{M} K \,\text{d}\mathcal{M} + \int_{\partial \mathcal{M}} \kappa_g \,\text{d}l & = & 0\\
    2 \pi \chi(\mathcal{M}) & = & 2\pi\\
  \end{array}
  \right.
\end{equation}
which is not in accordance with the Gauss-Bonnet theorem. Therefore, $J$ has to be singular somewhere in $\mathcal M$.

The goal is to build a usable parametrization of $\mathcal M$, \emph{i.e.}, being able to use this parametrization to build a quad mesh of $\mathcal M$. Therefore, we will allow $J$ to be singular on a finite number $N$ of points $\mathbf S_j$, $j\in[|0,N-1|]$ and show that this condition is sufficient for this problem to always have a unique solution.

\subsection{$H$ PDE at singular points}
\label{sec:localEquationAtSingularPoints}
For now, we know boundary conditions for $H$, Eq.~(\ref{eq:bndH}), and the local equation in smooth regions, Eq.~(\ref{eq:smoothH}). The only thing left is to determine a local PDE governing $H$ at singular points $\{\mathbf{S}_j\}$.
We define $k_j$ as the index of singularity $\mathbf{S}_j$.

\begin{figure}[h!t]
  \centering
  \includegraphics[width=0.25\textwidth]{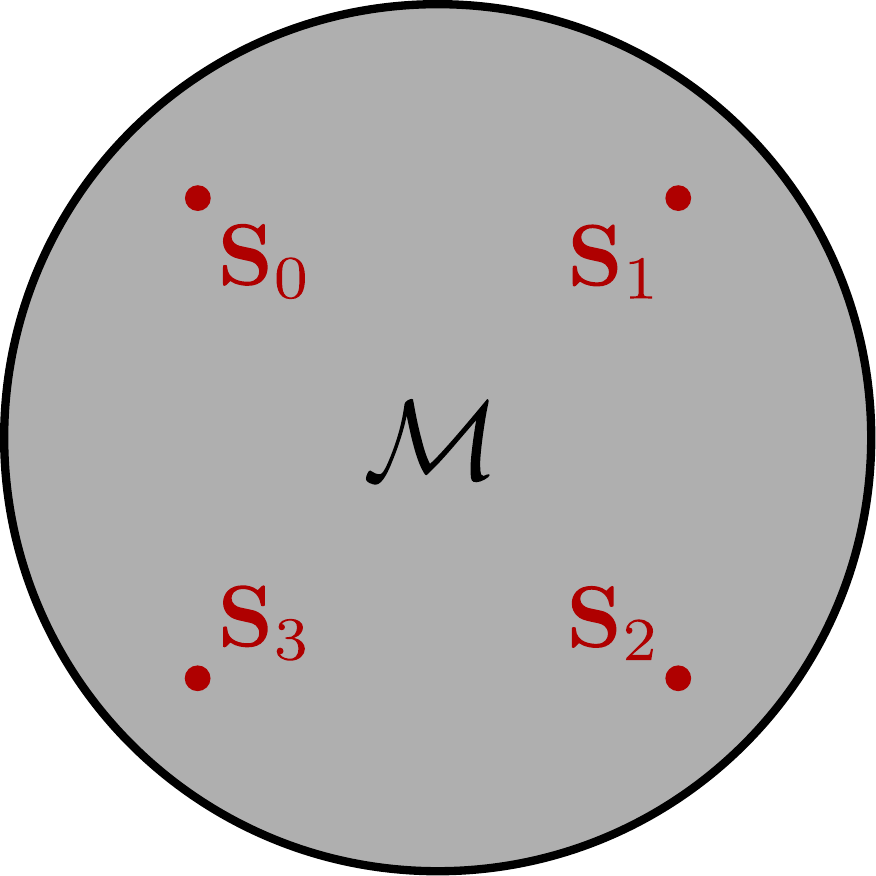}
  \caption{The disk with four singularities of index $1$.}
  \label{fig:4sing1}
\end{figure}

For this, we are making two reasonable assumptions:
\begin{equation}
  \left\{
  \begin{array}{rlll}
    \Delta H(\mathbf{S}_j) & = & -K(\mathbf{S}_j) + \alpha_j \delta (\mathbf{S}_j) & \\
%    \text{ where } \alpha_j \text{ is a constant} 
%     & & & \text{ and } \delta \text{ is Dirac distribution}\\
    k_i & = & k_j \Rightarrow \alpha_i=\alpha_j,&\\
  \end{array}
  \right.
  \label{eq:hypDirac}
\end{equation}

where $\alpha_j$ is a constant,  and  $\delta$ is the Dirac distribution. We consider the disk $\mathcal M$ represented in Fig.~\ref{fig:4sing1} with $4$ singularities $\mathbf{S}_j, j\in[|0,3|]$ of index $k_j=1$.

The Gauss-Bonnet theorem states that:
\begin{equation*}
  \int_\mathcal M K \,\text{d}\mathcal M + \int_{\partial \mathcal M} \kappa_g \text{d}l = 2 \pi \chi(\mathcal M).
\end{equation*}
Replacing $K$ and $\kappa_g$ by their values in Eq.~(\ref{eq:bndH}) and (\ref{eq:smoothH}), and using the hypothesis (\ref{eq:hypDirac}) we get $\alpha  =   2 \pi\frac{1}{4}$.

%\begin{figure*}[h!t]
%\begin{center}
%\begin{equation*}
%  \begin{array}{lrll}
%    & \int_\mathcal M -\Delta H\,\text{d}\mathcal M + \sum_i \alpha_i \delta (\mathbf{S}_i) \,\text{d}\mathcal M + \int_{\partial \mathcal M} -\nabla_{\mathbf t} H \text{d}l & = & 2 \pi \chi(\mathcal M)\\
%    \Rightarrow & \int_\mathcal M -\Delta H \,\text{d}\mathcal M - \int_{\partial \mathcal M} \nabla_{\mathbf t} H \text{d}l + \alpha\int_\mathcal M \sum_i \delta (\mathbf{S}_i) \,\text{d}\mathcal M & = &  2 \pi\\
%    \Rightarrow & \alpha & = &  2 \pi\frac{1}{4}\\
%  \end{array}
%\end{equation*}
%\end{center}
%\end{figure*}

For the singularity of index $1$ we have:
\begin{equation*}
    \Delta H(\mathbf{S}_j) = -K(\mathbf{S}_j) + 2\pi\frac{1}{4} \delta (\mathbf{S}_j).
\end{equation*}
Using the same idea, we can generalize the following:
\begin{equation}
    \Delta H(\mathbf{S}_j) = -K(\mathbf{S}_j) + 2\pi\frac{k_j}{4} \delta (\mathbf{S}_j).
\end{equation}

\subsection{Boundary value problem for $H$}
\label{sec:pdeH}
To sum up, the equations governing $H$ on $\mathcal M$ are:
\begin{equation}
  \left\{
  \begin{array}{llll}
    \Delta H & = & -K + 2\pi\frac{k_j}{4} \delta (\mathbf{S}_j) & \text{on } \mathcal M\\
    \nabla_\mathbf t H & = & -\kappa_g & \text{on } \partial\mathcal M.\\
  \end{array}
  \right.
  \label{eq:pdeH}
\end{equation}

\begin{figure}[h!t]
\begin{center}
  \includegraphics[width=0.4\textwidth]{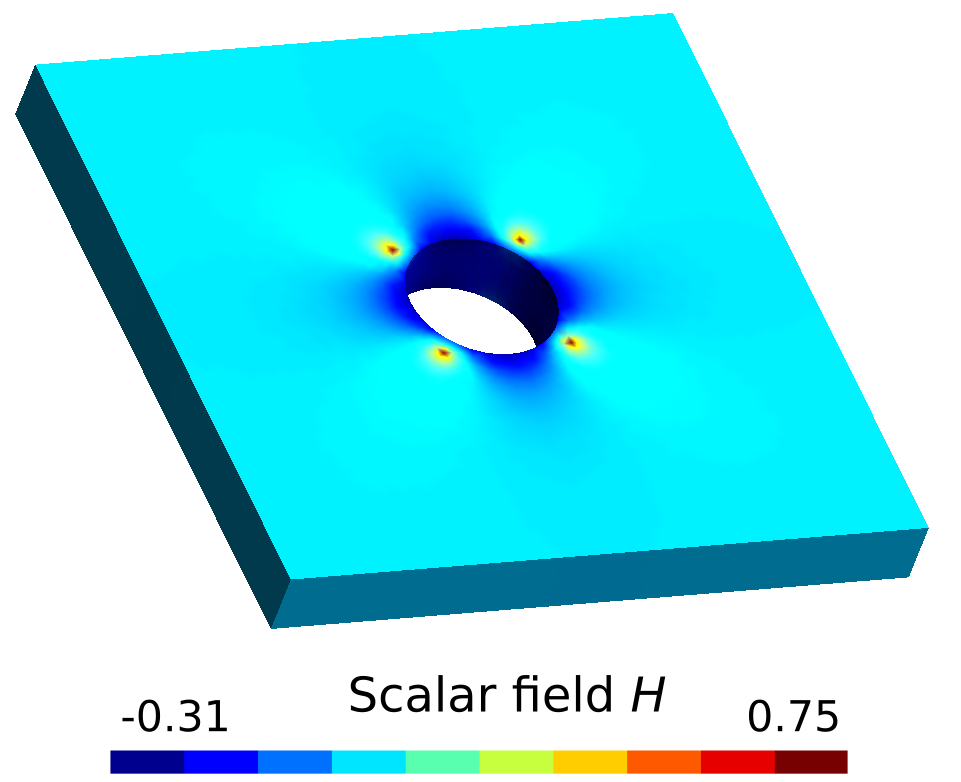}
\caption{$H$ function obtained on a closed manifold.}
\label{fig:hmanifold}
\end{center}
\end{figure}

This problem is well-posed and admits a unique solution to an arbitrary additive constant. A triangulation $\mathcal M_T$ of the manifold $\mathcal M$ is generated and problem~(\ref{eq:pdeH}) is solved using a finite element formulation with order 1 Lagrange elements. Once $H$ is determined (illustrated in Fig.~\ref{fig:hmanifold}), the next step is to retrieve $J$ orientation, detailed in the next section. The fact that $H$ is only known up to an additive constant is not harmful as only $\nabla H$ will be needed to retrieve $J$ orientation.

\subsection{Retrieving crosses orientation from H}
\label{sec:retrieveCross}
In order to get an orientation at a given point $\mathbf{X}\in\mathcal{M}$, a local reference basis $(\mathbf{u}_\mathbf{X},\mathbf{v}_\mathbf{X},\mathbf{n})$ in $\mathbf{X}$ is recalled.

Equation~(\ref{eq:intManifold2}) imposes that:
\begin{equation}
  \left\{
  \begin{array}{crcl}
    \nabla_{\mathbf{u}} H & = & \nabla_{\mathbf{v}}(\phi+\theta)\, + \cgamma \nabla_{\mathbf{v}}\psi\\
    \nabla_{\mathbf{v}} H & = & - \nabla_{\mathbf{u}}(\phi+\theta)\, - \cgamma \nabla_{\mathbf{u}}\psi\\
  \end{array}
  \right.
  \label{eq:intManifoldTheta1}
\end{equation}
which is equivalent to:
\begin{equation}
  \left\{
  \begin{array}{crcl}
    \nabla_{\mathbf{u}_\mathbf{X}} H & = & \nabla_{\mathbf{v}_\mathbf{X}}(\phi+\theta)\, + \cgamma \nabla_{\mathbf{v}_\mathbf{X}}\psi\\
    \nabla_{\mathbf{v}_\mathbf{X}} H & = & - \nabla_{\mathbf{u}_\mathbf{X}}(\phi+\theta)\, - \cgamma \nabla_{\mathbf{u}_\mathbf{X}}\psi\\
  \end{array}
  \right.
  \label{eq:intManifoldTheta2}
\end{equation}
and eventually gives:
\begin{equation}
  \left\{
  \begin{array}{crcl}
    \nabla_{\mathbf{u}_\mathbf{X}} \theta & = & - \nabla_{\mathbf{v}_\mathbf{X}} H - \nabla_{\mathbf{u}_\mathbf{X}}\phi\, - \cgamma \nabla_{\mathbf{u}_\mathbf{X}}\psi\\
    \nabla_{\mathbf{v}_\mathbf{X}}\theta & = & \nabla_{\mathbf{u}_\mathbf{X}} H - \nabla_{\mathbf{v}_\mathbf{X}}\phi\, - \cgamma \nabla_{\mathbf{v}_\mathbf{X}}\psi\\
  \end{array}
  \right.
  \label{eq:intManifoldTheta3}
\end{equation}
which is linear in $\theta$.

We can show that there always exists a scalar field $\theta$ verifying Eq.~(\ref{eq:intManifoldTheta3}). The $\theta$ exists if and only if we have:
\begin{equation}
  \nabla_{\mathbf{u}_\mathbf{X}}\nabla_{\mathbf{v}_\mathbf{X}} \theta - \nabla_{\mathbf{v}_\mathbf{X}}\nabla_{\mathbf{u}_\mathbf{X}} \theta = 0.
  \label{eq:rotThetaNul}
\end{equation}
Using Eq.~(\ref{eq:intManifoldTheta3}) we obtain:
\begin{equation}
  \begin{array}{ll}
    & \nabla_{\mathbf{u}_\mathbf{X}}\nabla_{\mathbf{v}_\mathbf{X}} \theta - \nabla_{\mathbf{v}_\mathbf{X}}\nabla_{\mathbf{u}_\mathbf{X}} \theta\\
    = &\Delta H + \nabla_{\mathbf{v}_\mathbf{X}}(c_\gamma \nabla_{\mathbf{u}_\mathbf{X}}\psi) - \nabla_{\mathbf{u}_\mathbf{X}}(c_\gamma\nabla_{\mathbf{v}_\mathbf{X}}\psi)\\
    = & -K + \nabla_{\mathbf{v}_\mathbf{X}}(c_\gamma \nabla_{\mathbf{u}_\mathbf{X}}\psi) - \nabla_{\mathbf{u}_\mathbf{X}}(c_\gamma\nabla_{\mathbf{v}_\mathbf{X}}\psi).\\
  \end{array}
  \label{eq:rotTheta}
\end{equation}
We know that, for 2D manifolds embedded in $\mathbb R^3$, the Gaussian curvature $K$ is equal to the jacobian of the Gauss map of the manifold \cite{singer2015lecture}. We have:
\begin{equation}
  \left\{
  \begin{array}{lll}
    \nabla_{\mathbf u _\mathbf X}\mathbf n & = & s_\gamma\nabla_{\mathbf u _\mathbf X}\psi\left(\begin{matrix}\cpsi\\\spsi\\0\end{matrix}\right)-\nabla_{\mathbf u _ \mathbf X}\gamma\left(\begin{matrix}-\spsi\cgamma\\\cpsi\cgamma\\\sgamma\end{matrix}\right)\\
        \nabla_{\mathbf v _\mathbf X}\mathbf n & = & s_\gamma\nabla_{\mathbf v _\mathbf X}\psi\left(\begin{matrix}\cpsi\\\spsi\\0\end{matrix}\right)-\nabla_{\mathbf v _ \mathbf X}\gamma\left(\begin{matrix}-\spsi\cgamma\\\cpsi\cgamma\\\sgamma\end{matrix}\right)\\
  \end{array}
  \right.
  \label{eq:jacGaussMap}
\end{equation}

Therefore we also have: 
\begin{equation}
  K = s_\gamma(\nabla_{\mathbf{v}_\mathbf{X}}\psi\nabla_{\mathbf{u}_\mathbf{X}}\gamma - \nabla_{\mathbf{u}_\mathbf{X}}\psi\nabla_{\mathbf{v}_\mathbf{X}}\gamma).
  \label{eq:gaussMap}
\end{equation}

Developing Eq.~(\ref{eq:rotTheta}) and substituting $K$ with the right-hand side of  Eq.~(\ref{eq:gaussMap}) we get:
\begin{equation}
  \begin{array}{ll}
    & -K + \nabla_{\mathbf{v}_\mathbf{X}}(c_\gamma \nabla_{\mathbf{u}_\mathbf{X}}\psi) - \nabla_{\mathbf{u}_\mathbf{X}}(c_\gamma\nabla_{\mathbf{v}_\mathbf{X}}\psi)\\
    = & -K + c_\gamma \nabla_{\mathbf{v}_\mathbf{X}}\nabla_{\mathbf{u}_\mathbf{X}}\psi - \sgamma\nabla_{\mathbf{v}_\mathbf{X}}\gamma \nabla_{\mathbf{u}_\mathbf{X}}\psi \\
    & - c_\gamma\nabla_{\mathbf{u}_\mathbf{X}}\nabla_{\mathbf{v}_\mathbf{X}}\psi + \sgamma\nabla_{\mathbf{u}_\mathbf{X}}\gamma\nabla_{\mathbf{v}_\mathbf{X}}\psi\\
    = & 0. \\
  \end{array}
  \label{eq:prrofRotTheta}
\end{equation}
As Eq.~(\ref{eq:rotThetaNul}) is verified, we know that there exists a scalar field $\theta$ verifying   Eq.~(\ref{eq:intManifoldTheta3}), and therefore that our problem has a unique solution.

In order to solve Eq.~(\ref{eq:intManifoldTheta3}), we first need to obtain a smooth global basis $(\mathbf{u}_\mathbf{X},\mathbf{v}_\mathbf{X},\mathbf{n})$ on $\mathcal{M}$. This is possible by generating a branch cut $\mathcal L$, as defined below, and computing a smooth global basis $(\mathbf{u}_\mathbf{X},\mathbf{v}_\mathbf{X},\mathbf{n})$ on $\mathcal{M}$ allowing discontinuities across $\mathcal L$.

A branch cut is a set $\mathcal{L}$ of curves of a domain $M$ that do not form any closed loop 
and that cut the domain in such a way that it is impossible to find any closed loop in $M\setminus \mathcal{L}$ that encloses one or several singularities, or an internal boundary. 
As we already have a triangulation of $M$, the branch cut $\mathcal{L}$ 
is in practice simply a set of edges of the triangulation.

\begin{figure}[h!t]
\begin{center}
  \includegraphics[width=0.4\textwidth]{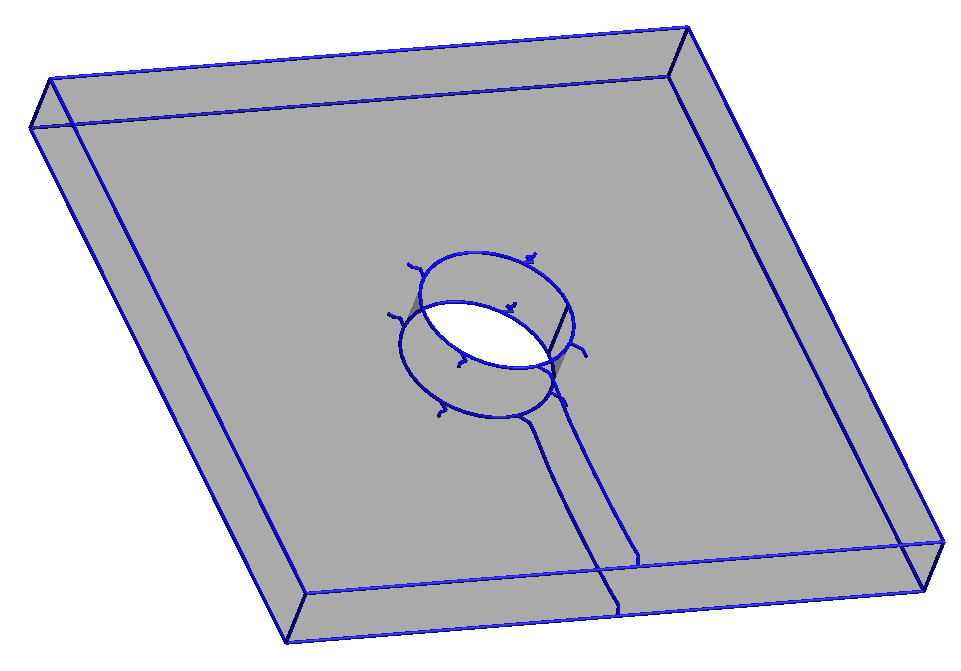}
\caption{Edges of the branch cut $\mathcal{L}$ are represented in blue. There exists no closed loop in $M \setminus \mathcal L$ enclosing one or several singularities.}
\label{fig:cutGraph}
\end{center}
\end{figure}

The branch cut is generated with the method described in \cite{Jezdimirovic:2021} which is based on \cite{bommes2009}. An example of generated branch cut is presented in Fig.~\ref{fig:cutGraph}.

Once a branch cut $\mathcal{L}$ is available, the field $\theta$ can be computed by solving the linear equations~(\ref{eq:intManifoldTheta3}). With equations~(\ref{eq:intManifoldTheta3}), $\theta$ is known up to an additive constant. For the problem to be well-posed, $\theta$ value has to be imposed at one point of domain $\mathcal M$. The chosen boundary condition consists in fixing the angle $\theta$ at {\em one} arbitrary point $\mathbf X_{BC} \in \partial \mathcal M$ so that  $\mathcal{C}_\mathcal{M} (\mathbf X_{BC})$ has one of its branches collinear with $\mathbf T (\mathbf X)$.
The problem can be rewritten as the well-posed Eq.~(\ref{eq:thetaManifold}) and is solved using the finite element method on the triangulation $\mathcal{M}_T$ with order one Crouzeix-Raviart elements. This kind of elements has shown to be more efficient for cross-field representation \cite{jezdimirovic2019}.
\begin{equation}
  \left\{\begin{array}{rl}
  P_{T\mathcal M}(\nabla \theta) & = P_{T\mathcal M}(\mathbf n \times \nabla H - \nabla \phi - \cgamma \nabla \psi) \text{ in } \mathcal M \\
  \theta(\mathbf X_{BC}) & = \theta_{\mathbf X_{BC}} \text{ for an arbitrary } \mathbf X_{BC} \in\partial \mathcal M\\
  \theta & \text{ discontinuous on } \mathcal{L}
\end{array}\right.
\label{eq:thetaManifold}
\end{equation}

It is important to note that for Eq.~(\ref{eq:thetaManifold}) to be well-posed, the $\theta$ value can only be imposed on a \emph{single} point. A consequence is that if $\mathcal M$ has more than one boundary ($\partial \mathcal M = \partial \mathcal M_1 \cup \partial\mathcal M_2 \cup \cdots \cup \partial\mathcal M_n$), the resulting cross-field is guaranteed to be tangent to the boundary $\partial\mathcal M_i$ such as $\mathbf X_{BC}\in \partial\mathcal M_i$, which does not necessarily hold for all boundaries $\partial M_j$ for $j\neq i$, as detailed in Section~\ref{sec:nonQuadMeshableSingConfig}.

\begin{figure}[h!t]
\begin{center}
  \includegraphics[width=0.4\textwidth]{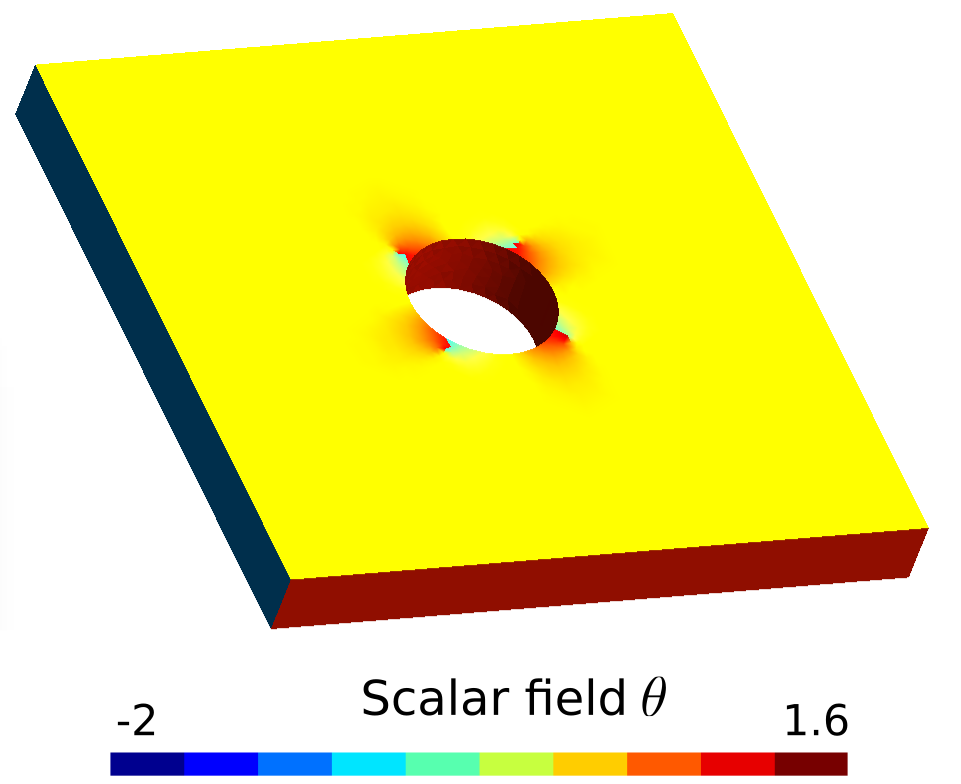}
\caption{Scalar field $\theta$ obtained from scalar field $H$ (represented in Fig.~\ref{fig:hmanifold}).}
\label{fig:thetamanifold}
\end{center}
\end{figure}

Once $H$  and $\theta$ scalar fields are computed on $\mathcal{M}$ (illustrated respectively in Fig.~\ref{fig:hmanifold} and Fig.~\ref{fig:thetamanifold}), the cross-field $\mathcal{C}_\mathcal{M}$ can be retrieved for all $\mathbf X \in \mathcal{M}$:
\begin{equation}
\mathbf c (\mathbf X) = \{ \mathbf u_k = \mathcal{R}_{\theta + k\frac{\pi}{2},\mathbf{n}}(\mathbf{u}_{\mathbf{X}}), k\in[|0,3|]\}.
\label{eq:HThetaToCross}
\end{equation}

\section{Preliminary results}
\label{S:Preliminary}

\begin{figure*}[h!t]
\begin{center}
\includegraphics[width=0.85\textwidth]{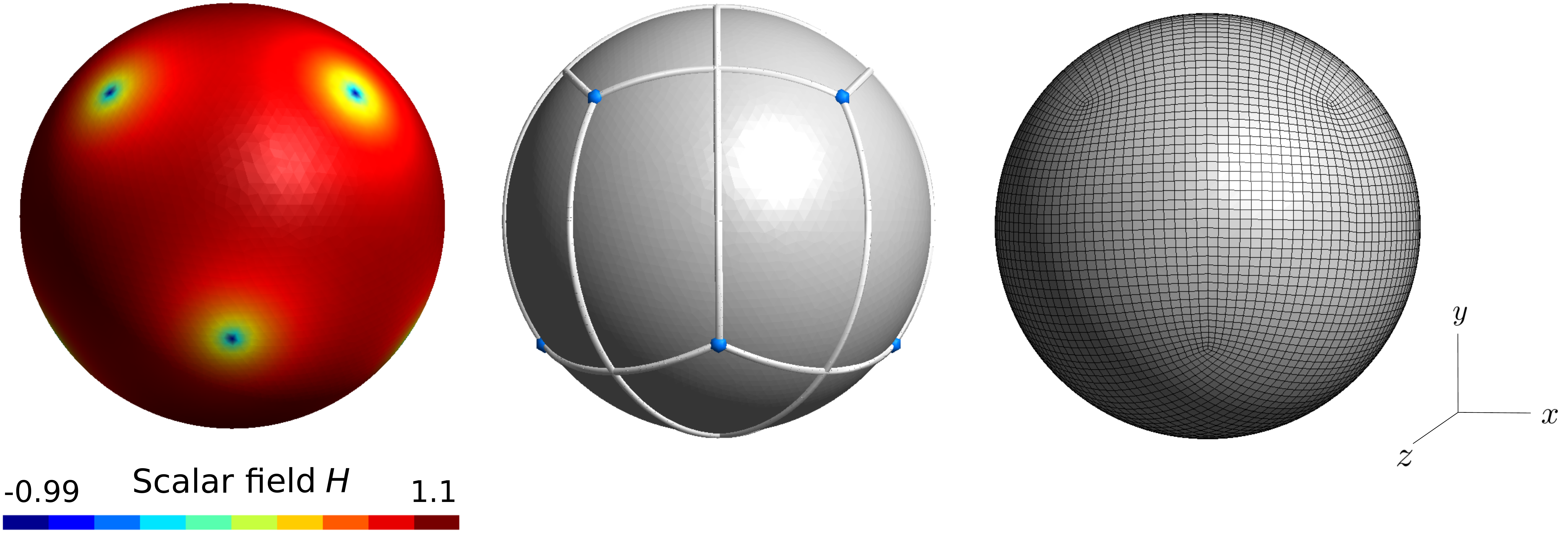}
\caption{Quad mesh on a 2-sphere with a natural singularity configuration forming an anticube. The singularity configuration comes from solving a non-linear problem, \emph{i.e.}, by using the MBO algorithm from \cite{viertel2019approach}.}
\label{fig:pipeline1}
\end{center}
\end{figure*}

\begin{figure*}[h!t]
\begin{center}
\includegraphics[width=0.85\textwidth]{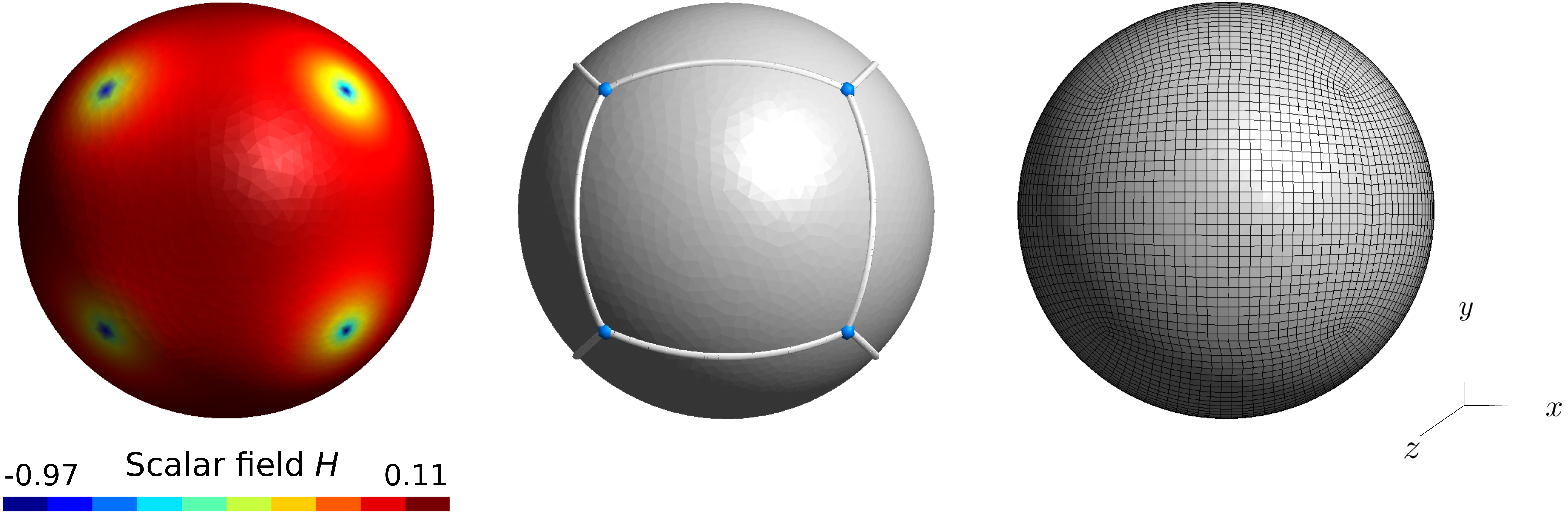}
\caption{Quad mesh on a 2-sphere with an imposed singularity configuration forming a cube.}
\label{fig:pipeline2}
\end{center}
\end{figure*}

As a proof of concept, the cross-field computation based on imposed singularity configuration is included in the 3-step quad meshing pipeline of \cite{Jezdimirovic:2021} (illustrated in Fig.~\ref{fig:pipeline1} and Fig.~\ref{fig:pipeline2}): 

\noindent{\textbf{Step 1:}} impose a singularity configuration, \emph{i.e.}, position and valences of singularities (see \cite{Jezdimirovic:2021}).

\noindent{\textbf{Step 2:}} compute a cross-field with the prescribed singularity configuration of Step 1 on an adapted mesh (singularities are placed in refined regions), by solving only two linear systems (Section~\ref{S:Isotropic}).

\noindent{\textbf{Step 3:}} compute a quad layout on the accurate cross-field of Step 2, and generate a full block-structured isotropic quad mesh (see \cite{jezdimirovic2019, Jezdimirovic:2021}).

The presented pipeline includes the automatic check that singularity configuration obeys the Euler characteristic of the surface, but it does not inspect all Abel-Jacobi conditions \cite{Chen:2019, Lei:2020, zheng2021quadrilateral}. 
Further, the models of industrial complexity would require a more robust quad layout generation technique than the one followed here (\cite{jezdimirovic2019, Jezdimirovic:2021}).
The final quad mesh is isotropic, obtained from the quad layout via per-partition bijective parameterization aligned with the smooth cross-field (singularities can only be located on corners of the partitions) \cite{Jezdimirovic:2021}, and following the size map implied by the $H$, \emph{i.e.}, the element's edge length is $s=e^H$.
In case when the application demands an anisotropic quad mesh, two sizing fields $(H_1, H_2)$ for the cross-field must be computed, more details in Section~\ref{S:Anisotropic}. 


\subsection{Valid singularity configurations for quad meshing} 
\label{sec:nonQuadMeshableSingConfig}

\begin{figure}[h!t]
\begin{center}
  \includegraphics[width=0.5\textwidth]{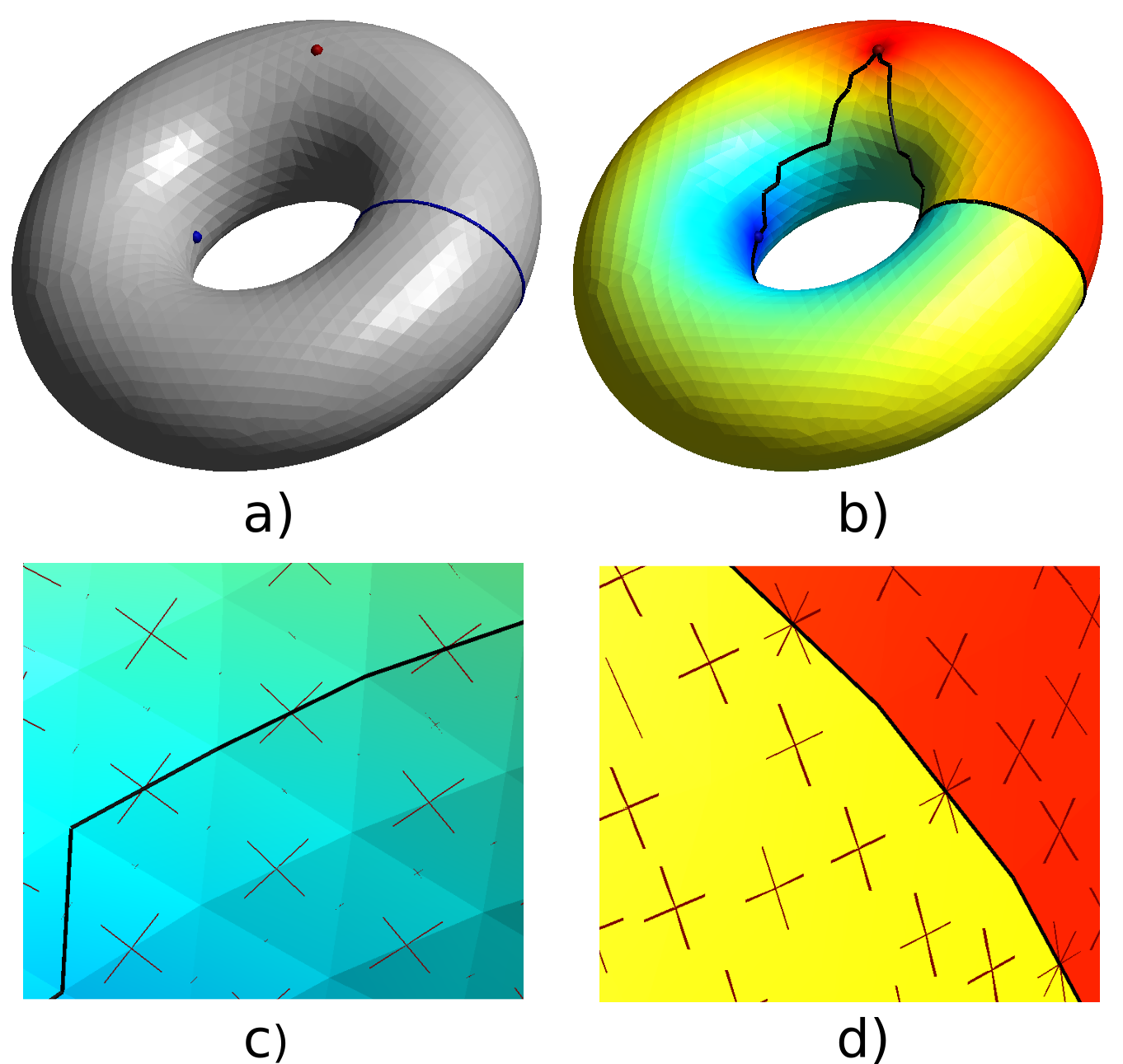}
\caption{Imposing a $3-5$ singularity configuration on a torus with the boundary marked in blue a), marking the boundary and the cut graph in black b), consistent cross-field across the cut graph c), and cross-field not aligned with the boundary d).}
\label{fig:AJplanar1}
\end{center}
\end{figure}

The singularity configuration, including both the positions and valences, plays a crucial role in the generation of conformal quad meshes \cite{gu2020computational}. It is essential to note that not all user-imposed singularity configurations matching the Euler's characteristic of the surface will be valid for quad meshing, Fig.~\ref{fig:AJplanar1}. The central cause for this lies in the fact that a combination of choices of valences and holonomy is not arbitrary \cite{mylesrobust}. Relevant findings on the non-existence of certain quadrangulations can be found in \cite{Barnette:1971, Jucovic:1973, Izmestiev:2013}.

The work of \cite{beaufort2017computing} presents the formula for determining the numbers of and indices of singularities, and \cite{fogg2018singularities} their possible  combinations in conforming quad meshes. Latter authors also show that the presented formula is necessary but not sufficient for quad meshes, but neither of these works are proving the rules for the singularities’ placement.

Recently, the sufficient and necessary conditions for valid singularity configuration of the quad mesh are presented in the framework based on \textit{Abel-Jacobi's theory} \cite{Chen:2019, Lei:2020, zheng2021quadrilateral}. The developed formulation here is under its direct constraint. In practice, imposing a singularity configuration fulfilling Euler's characteristic constraint ensures that the flat metric, \emph{i.e.}, the $H$ field can be obtained. If this singularity configuration also verifies the holonomy condition, the cross-field will be aligned with all boundaries and consistent across the cut graph.

We recall here that our formulation entitles the user to impose its own singularity configuration, which in practice can contain a suboptimal distribution of singularities. As a consequence, computed cross-field may not be aligned with all boundaries, disabling the generation of the final conformal isotropic quad mesh. To bypass this issue, the following section develops an integrable cross-field formulation with two independent metrics (which are flat except at singularities), instead of only one as presented for Abel-Jacobi conditions.

\subsection{Dealing with suboptimal distribution of singularities} 
\label{sec:SubOptimalSing}

\begin{figure*}[h!t] 
  \begin{center}
    \includegraphics[width=0.9\textwidth]{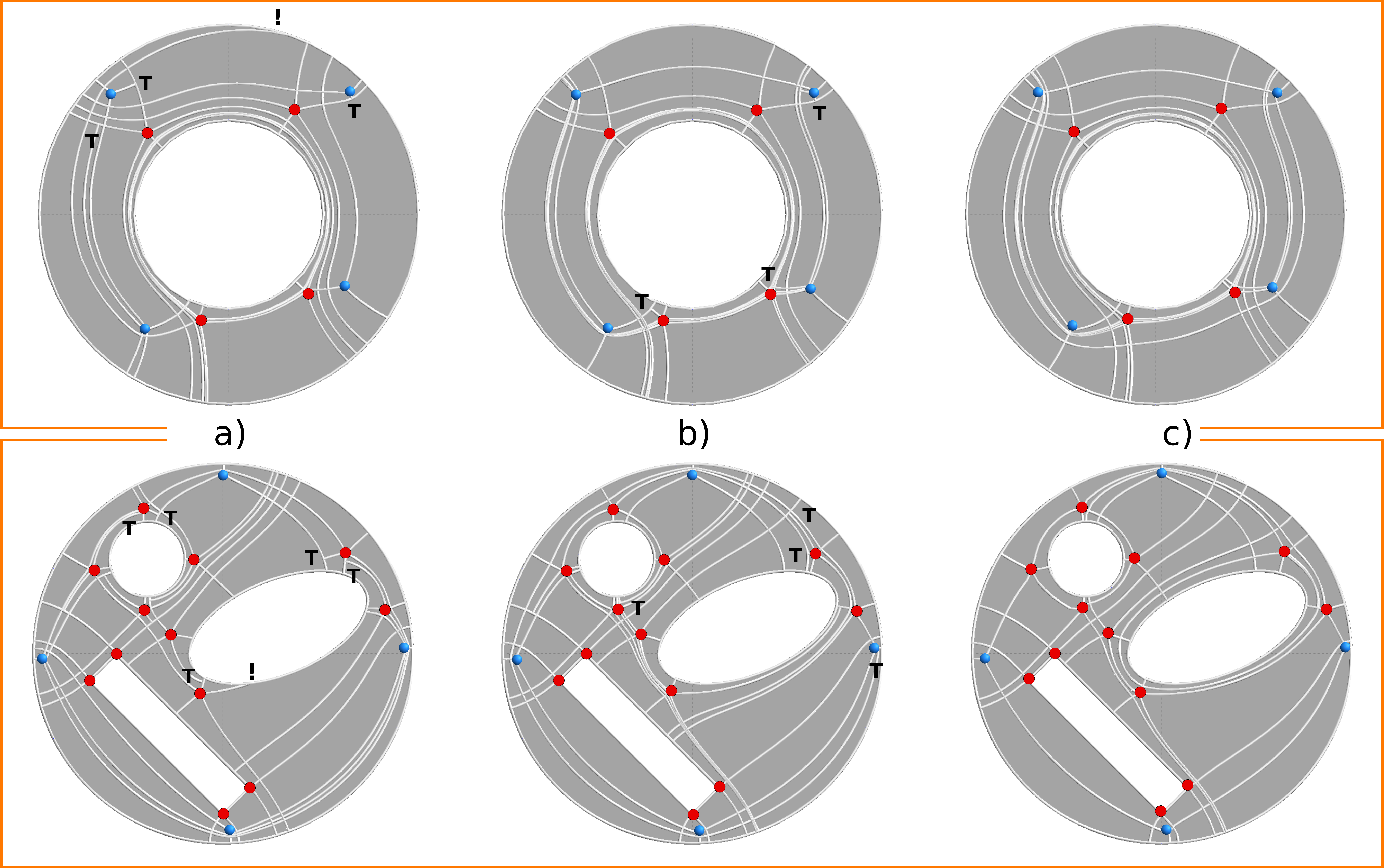}
    \caption{Obtained quad layouts on an imposed set of singularities that do not respect the location’s condition from the Abel-Jacobi theorem. Quad layouts obtained using the integrable cross-field with isotropic scaling are not aligned with boundaries (marked with "$\textbf{!}$") and  demonstrate the presence of t-junctions (marked with "\textbf{T}") a). Quad layouts obtained with imposing the $\theta$ value are boundary aligned but not integrable and demonstrate the presence of t-junctions b). Quad layouts obtained using the integrable cross-field with an anisotropic scaling c).}
    \label{fig:prDef} 
  \end{center}
\end{figure*}

The issue of suboptimal distribution of singularities imposes the need for developing a new cross-field formulation on the imposed singularity configuration, which considers the integrability while relaxing the condition on isotropic scaling of crosses' branches. More specifically, the integrability condition, along with computing only one scaling field $H$, $||\tilde{\mathbf u}||=||\tilde{\mathbf v}||$, imposes the strict constraint on the valid singularity configurations, \emph{i.e.}, the need for fulfilling the \textit{Abel-Jacobi theorem}. Therefore, two sizing fields $L_1 = ||\tilde{\mathbf u}||$ and $L_2=||\tilde{\mathbf v}||$ are introduced and the upcoming section presents the mathematical foundations for the generation of an integrable cross-field with anisotropic scaling on $2-$D manifolds. As it will be shown in the following, this setting presents promising results in generating an integrable and boundary-aligned cross-field on the imposed set of singularities, even when their distribution is not fulfilling all Abel-Jacobi conditions. Only for the sake of visual comprehensiveness, the presented motivational examples in Fig.~\ref{fig:prDef} - Fig.~\ref{fig:ex2IntConv} are planar.

\section{Integrability condition with anisotropic scaling}
\label{S:Anisotropic}

As explained previously, a cross-field $\mathcal C_{\mathcal{M}}$ is integrable if and only if $\tilde{\mathbf u}$ and $\tilde{\mathbf v}$ commute under the Lie Bracket. In other words, the condition:
\begin{equation}
0 = [\tilde{\mathbf{u}},\tilde{\mathbf{v}}] 
= \nabla_{\tilde{\mathbf u}}  \tilde{\mathbf v} -\nabla_{\tilde{\mathbf v}} \tilde{\mathbf u} 
= [ L_1 \mathbf u , L_2 \mathbf v] 
\label{eq:chap3:integrability}
\end{equation}
where:
\begin{equation}
  L_1 = ||\tilde{\mathbf u}||\text{, } \quad \quad L_2 = ||\tilde{\mathbf v}||
\end{equation}
and
\begin{equation}
  \begin{array}{rcl}
    \mathbf u & = & \frac{\tilde{\mathbf u}}{||\tilde{\mathbf u}||}\\
    \mathbf v & = & \frac{\tilde{\mathbf v}}{||\tilde{\mathbf v}||}\\
  \end{array}
\end{equation}
has to be verified.

Developing the latter expression and posing for convenience $L_1 = e^{H_1}$ and $L_2 = e^{H_2}$, it becomes: 
$$
0 = \mathbf v \nabla_{\mathbf u} H_2 - \mathbf u \nabla_{\mathbf v} H_1 + [\mathbf u , \mathbf v],
$$
and then
\begin{equation}
  \left\{
  \begin{array}{crcl}
    \nabla_\mathbf{u} H_2 & = & - < \mathbf v,  [\mathbf u , \mathbf v]  > \\
    \nabla_\mathbf{v} H_1 & = & < \mathbf u,  [\mathbf u , \mathbf v]  >
  \end{array}
  \right.
  \label{eq:chap3:intManifold}
\end{equation}
which after the substitution of Eq.~(\ref{bracketuv}) gives:
\begin{equation}
  \left\{
  \begin{array}{rrl}
    \nabla_{\mathbf u} H_2 & = &  \nabla_{\mathbf v}\theta\, + \nabla_{\mathbf v}\phi\, + \cgamma \nabla_{\mathbf v}\psi \\
    -\nabla_{\mathbf v} H_1 & = & \nabla_{\mathbf u}\theta\, + \nabla_{\mathbf u}\phi\, + \cgamma \nabla_{\mathbf u}\psi.
  \end{array}
  \right.
  \label{eq:chap3:intManifold2}
\end{equation}

It is important to note that the three scalar fields $(\theta,H_1,H_2)$ are completely defining the cross-field $\mathcal C_{\mathcal{M}}$, as $(\psi,\gamma,\phi)$ are known since they are defining the local manifold basis $(\mathbf t, \mathbf T, \mathbf n)$.

From Eq.~(\ref{eq:chap3:intManifold2}), we can define the cross-field $\mathcal C_{\mathcal{M}}$ integrability error $E$ as:
\begin{equation}
\begin{array}{ll}
 & E^2(\theta,H_1,H_2)  \\
 = & \int_{\mathcal{M}} (\nabla_{\vec u} H_2 - \nabla_{\vec v}\theta\, - \nabla_{\vec v}\phi\, - \cgamma \nabla_{\vec v}\psi)^2 \\
& \, \, + (\nabla_{\vec v} H_1 + \nabla_{\vec u}\theta\, + \nabla_{\vec u}\phi\, + \cgamma \nabla_{\vec u}\psi)^2\,\text{d}{\mathcal{M}}.  
  \end{array}
  \label{eq:intError}
\end{equation}

The problem of generating an integrable cross-field with anisotropic scaling can therefore be reduced at finding three scalar fields $(\theta,H_1,H_2)$ verifying $E(\theta,H_1,H_2)=0$.

The process of solving this problem presents several difficulties. First, the quadruple $(\theta,\psi,\gamma,\phi)$ are multivalued functions. This kind of difficulty is commonly encountered in cross-field generation and is tackled here by cutting the domain $\mathcal{M}$ along a generated cut graph. Then, minimizing $E$ regarding $(\theta,H_1,H_2)$ is an ill-posed problem. Indeed, there are no constraints on $\nabla_{\mathbf u} H_1$ and $\nabla_{\mathbf v} H_2$. This is the main obstacle for generating an integrable 2D cross-field with an anisotropic scaling.

A simple approach to solve this problem is proposed here. In order to do so, it is needed to:
\begin{itemize}
\setlength\itemsep{-0.5em}
\item be able to generate a boundary-aligned cross-field matching the imposed singularity configuration, 
\item compute $(H_1,H_2)$ minimizing $E$ for an imposed $\bar \theta$, 
\item compute $\theta$ minimizing $E$ for an imposed $(\bar H_1,\bar H_2).$ 
\end{itemize}
 
The final resolution solver (Algorithm~\ref{alg:generalH1H2theta}), proposed in Section~\ref{sec:minGeneralInt}, allows finding a local minimum for $E$ around an initialization $(\theta^0,H_1^0,H_2^0)$.

\subsection{Local manifold basis generation and $\theta$ initialization}
\label{sec:thetaInit}

As exposed earlier, in order to completely define a unitary cross-field $\mathcal C_{\mathcal{M}}$ with a scalar field $\theta$ it is needed to define a smooth global basis $(\mathbf t,\mathbf T,\mathbf n)$ on $\mathcal{M}$. This is possible by generating a branch cut $\mathcal L$ and computing a smooth global basis $(\mathbf t,\mathbf T,\mathbf n)$ on $\mathcal{M}$ allowing discontinuities across $\mathcal L$.

The branch cut is generated using the method described in \cite{bommes2009}. A local basis $(\mathbf t,\mathbf T,\mathbf n)$ on $\mathcal{M}$ can be generated with any cross-field method. Such local basis will be smooth and will not show any singularities, as discontinuities are allowed across the cut graph $\mathcal L$ and no boundary alignment  is required. Once the cut graph $\mathcal L$ and the local basis $(\mathbf t,\mathbf T,\mathbf n)$ are generated, it is possible to compute $\theta$ only if:
\begin{itemize}
\setlength\itemsep{-0.5em}
\item $\theta$ values on $\partial \mathcal{M}$ are known,
\item $\theta$ jump values across $\mathcal L$ are known.
\end{itemize}

These can be found using methods described in \cite{bommes2009}, or can be deduced from a low computational cost cross-field generation detailed in \cite{Jezdimirovic:2021}.

\subsection{Computing $(H_1,H_2)$ from imposed $\bar \theta$}
For a given $\bar \theta$, it is possible to find $(H_1,H_2)$ minimizing $E$. It is important to note that, in general, there does not exist a couple $(H_1,H_2)$ such as $E=0$. Minimizing $E$ with imposed $\bar \theta$ is finding the couple $(H_1,H_2)$ for which the integrability error is minimal.

The problem to solve is the following:               
\begin{equation}  
\begin{array}{ll} 
\text{Find }(\bar H_1,\bar H_2)\text{ such as } & \\                                                                                       
E(\bar \theta,\bar H_1,\bar H_2)=\displaystyle \min_{(H_1,H_2)\in(\mathcal C^1(\mathcal{M}))^2} E(\bar \theta,H_1,H_2). &   
\label{eq:pbIntH1H2theta} 
\end{array}                                                                            
\end{equation}                                                                                          
Let's define $\mathcal S$ as:
 $$\mathcal S = \{(\bar H_1, \bar H_2) \mid  (\bar H_1, \bar H_2) \text{ verifies Eq. (} \ref{eq:pbIntH1H2theta} \text{)} \}.$$                                                                                      
                                                                                                        
For this problem to be well-posed, a necessary condition is to have $2$ independent scalar equations involving $\nabla H_1$, and the same for $\nabla H_2$. We can note that in our case, there are no constraints on $\nabla_{\vec u} H_1$ and $\nabla_{\vec v} H_2$. Therefore, there is only $1$ scalar equation involving $\nabla H_1$, and $1$ scalar equation involving $\nabla H_2$. As a consequence, the problem we are looking to solve is ill-defined. As this problem is ill-defined, $\mathcal S$ will not be a singleton and, in the general case, there will be more than one solution to the problem (\ref{eq:pbIntH1H2theta}).

\begin{figure}[h] 
  \begin{center}
    \includegraphics[width=0.25\textwidth]{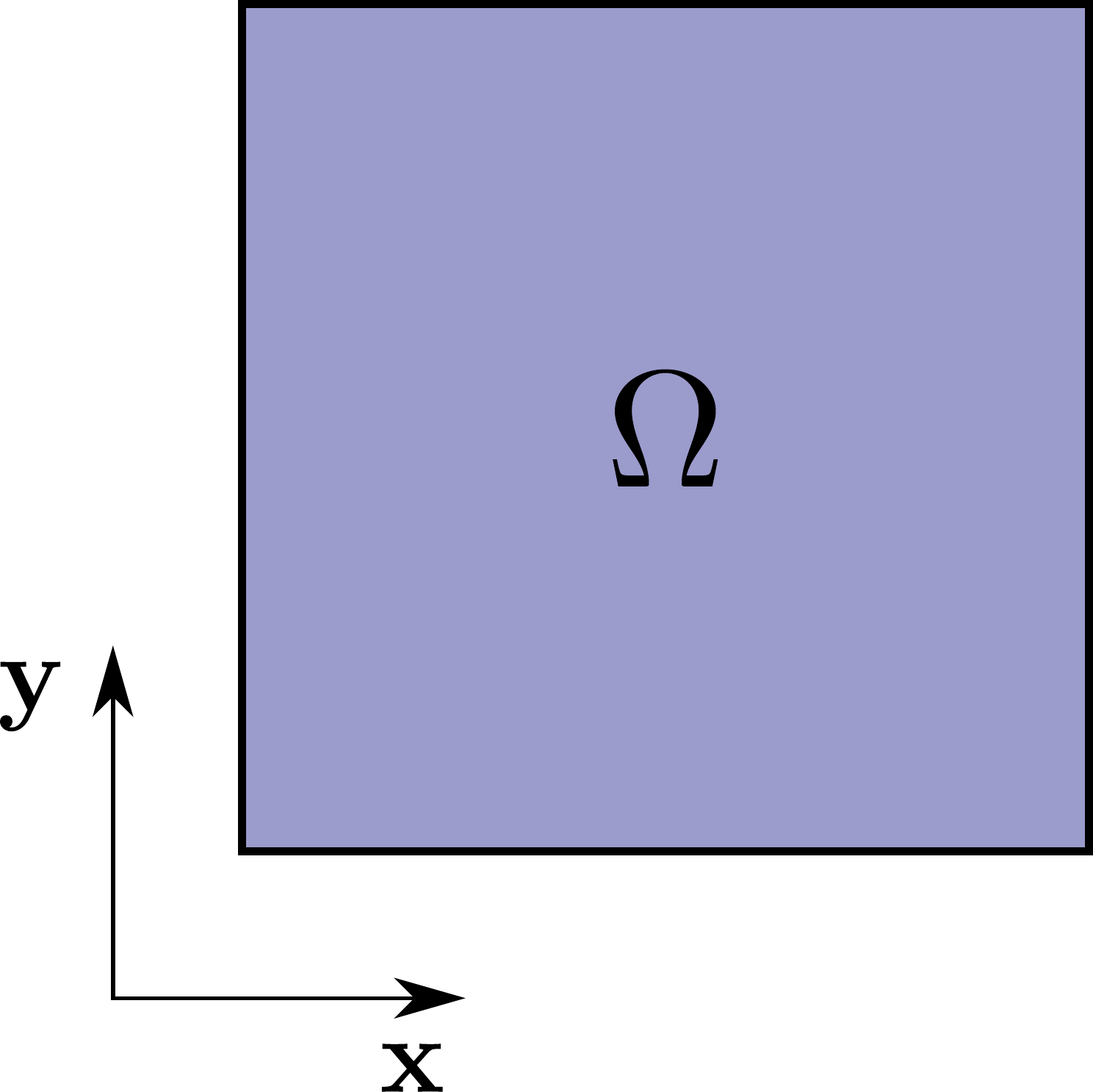}
    \caption{Planar square.}
    \label{fig:square} 
  \end{center}
\end{figure}

To discuss this problem in detail, we will use the simple example of a planar domain $\Omega$ illustrated in Fig.~\ref{fig:square}.

In this case, the unitary frame field $\mathcal C_{\Omega}$ obtained with common methods is:
\begin{equation}
  \mathcal{C}_{\Omega} = \{c(\vec X)=\{ \vec x, \vec y, -\vec x, -\vec y\}, \vec X\in{\Omega}\}
\end{equation}
which is equivalent to:
\begin{equation}
\bar\theta = 0.
\end{equation}
As in this case domain $\Omega$ is planar, we also have:
\begin{equation}
\psi = \gamma = \phi = 0.
\end{equation}

Equation~(\ref{eq:chap3:intManifold2}) becomes:                                                           
                                                                                                        
\begin{equation}                                                                                        
  \left\{                                                                                               
  \begin{array}{rrl}                                                                                    
    \nabla_{\vec x} H_2 & = &0 \\                                                                       
    -\nabla_{\vec y} H_1 & = & 0                                                                        
  \end{array}                                                                                           
  \right.                                                                                               
  \label{eq:chap3:intExSquare}                                                                          
\end{equation}                                                                                          
which gives:                                                                                            
\begin{equation}                                                                                        
  \left\{                                                                                               
  \begin{array}{rrl}                                                                                    
    H_1(x, y) & = & f(x)\text{, } \forall (x, y)\in\Omega\text{, }\forall f\in\mathcal{C}^1(\mathbb{R})\\
    H_2(x, y) & = & g(y)\text{, } \forall (x, y)\in\Omega\text{, }\forall g\in\mathcal{C}^1(\mathbb{R}). \\                                                                                                      
  \end{array}                                                                                           
  \right.                                                                                               
  \label{eq:chap3:intExSquare2}                                                                         
\end{equation}                                                                                          
                                                                                                        
Knowing this, we finally have $\mathcal S = \mathcal (C^1(\mathbb{R}))^2$. There is an infinity of solutions, confirming the fact that problem~(\ref{eq:pbIntH1H2theta}) is ill-defined.                          
                                                                                                        
The solution we could expect to obtain for quad meshing purposes would be:                              
\begin{equation}                                                                                        
  \mathcal S = \{ (H_1,H_2) = (0,0)\},                                                                  
  \label{eq:solSquare}                                                                                  
\end{equation}                                                                                          
which is equivalent to $(L_1,L_2)=(1,1)$.

Based on this simple example, we can deduce that problem (\ref{eq:pbIntH1H2theta}) has to be regularized in order to reduce the solution space. One way to achieve this goal is to add a constraint on the $(H_1,H_2)$ fields we are looking for. A natural one is to look for $(H_1,H_2)$ verifying Eq.~(\ref{eq:pbIntH1H2theta}) and being as smooth as possible.

With this constraint, the problem to solve becomes:
\begin{equation}
\begin{array}{ll}
 \text{Find }(\bar H_1,\bar H_2)\in\mathcal S\text{ such as } &\\
 \displaystyle\int_{\mathcal{M}} ||\nabla \bar H_1||^2 + ||\nabla \bar H_2||^2\,\text{d}\mathcal{M}  &\\ 
=  \displaystyle\min_{(H_1,H_2)\in\mathcal S} \displaystyle\int_{\mathcal{M}} ||\nabla H_1||^2 + ||\nabla H_2||^2\,\text{d}{\mathcal{M}}. & \\
  \end{array}
  \label{eq:pbIntH1H2thetaReg}
\end{equation}

Adding this constraint transforms the linear problem (\ref{eq:pbIntH1H2theta}) into a non-linear one (\ref{eq:pbIntH1H2thetaReg}). Algorithm~\ref{alg:H1H2reg} is used to solve Eq.~(\ref{eq:pbIntH1H2thetaReg}), leading to an $E$'s local minimum $(\bar \theta, \bar H_1, \bar H_2)$ close to $(\bar\theta, H_1^0, H_2^0)$.

\begin{algorithm}[h]
  \SetAlgoLined
  $k=0$\\
  initial guess $H_1^0$, $H_2^0$\\
  compute $\epsilon^0=E(\bar\theta,H_1^0,H_2^0)$\\
  \While{$\epsilon^k < \epsilon^{k-1}$}{
    $k=k+1$\\
    find $(H_1^k,H_2^k)$ minimizing:
    \begin{equation*}
    \begin{array}{rl}
     E(\bar\theta,f_1,f_2)+\int_{\mathcal{M}} & ||\nabla f_1 - \nabla H_1^{k-1}||^2 + \\
    & ||\nabla f_2 - \nabla H_2^{k-1}||^2\text{d}\,\mathcal{M}, \\
     \end{array}
      \end{equation*}
     $(f_1,f_2)\in\left(\mathcal{C}^1(\mathcal{M})\right)^2$\\
    compute $\epsilon^k=E(\bar\theta,H_1^k,H_2^k)$
  }
  \caption{Regularized solver for $(H_1,H_2)$}
  \label{alg:H1H2reg}
\end{algorithm}

\subsection{Computing $\theta$ from $(\bar H_1 \bar H_2)$}

For an imposed couple $(\bar H_1 \bar H_2)$, it is possible to find $\theta$ minimizing $E$. The problem to solve is formalized as:
\begin{equation}
\begin{array}{ll}
  \text{Find }\bar \theta \in \mathcal C^1(\mathcal{M})\text{ such as } & \\
   E(\bar \theta,\bar H_1,\bar H_2)=\displaystyle \min_{\theta\in\mathcal C^1(\mathcal{M})} E(\theta,\bar H_1,\bar H_2).&
  \label{eq:pbIntthetaH1H2}
\end{array}
\end{equation}

This problem is non-linear too since $\nabla_{\mathbf v} H_1$ and $\nabla_{\mathbf u} H_2$ are showing a non-linear dependence regarding $\theta$. Algorithm~\ref{alg:theta} is used to solve Eq.~(\ref{eq:pbIntthetaH1H2}), leading to an $E$'s local minimum $(\bar \theta, \bar H_1, \bar H_2)$ close to $(\theta^0, \bar H_1, \bar H_2)$.

\begin{algorithm}[h]
  \SetAlgoLined
  $k=0$\\
  initial guess $\theta^0$\\
  deduce $(\mathbf u^0,\mathbf v^0)$ from $\theta^0$\\
  compute $\epsilon^0=E(\theta^0,\bar H_1,\bar H_2)$\\
  \While{$\epsilon^k < \epsilon^{k-1}$}{
    $k=k+1$\\
    find $\theta^k$ minimizing:
    \begin{equation*}
    \begin{array}{rl}
   & E^k(f,\bar H_1,\bar H_2) \\
    = & \int_{\mathcal{M}} (\nabla_{\mathbf u^{k-1}} \bar H_2 - \nabla_{\mathbf v^{k-1}}f\, - \nabla_{\mathbf v^{k-1}}\phi\, - \cgamma \nabla_{\mathbf v^{k-1}}\psi)^2+ \\
 & \, \quad  (\nabla_{\mathbf v^{k-1}} \bar H_1 + \nabla_{\mathbf u^{k-1}}f\, + \nabla_{\mathbf u^{k-1}}\phi\, + \cgamma \nabla_{\mathbf u^{k-1}}\psi)^2\,\text{d}\mathcal{M}\\
%& f\in\mathcal{C}^1(\mathcal{M})  
     \end{array}
      \end{equation*}
    $f\in\mathcal{C}^1(\mathcal{M})$ \\
    deduce $(\mathbf u^k,\mathbf v^k)$ from $\theta^k$\\
    compute $\epsilon^k=E(\theta^k,\bar H_1,\bar H_2)$\;
  }
  \caption{Solver for $\theta$}
  \label{alg:theta}
\end{algorithm}

\subsection{Minimizing integrability error $E$ regarding $(\theta,H_1,H_2)$}
\label{sec:minGeneralInt}
Using the three steps exposed previously, it is possible to find a local minimum in the vicinity of an initialization $(\theta^0,H_1^0,H_2^0)$ following Algorithm~\ref{alg:generalH1H2theta}.

\begin{algorithm}[h]
\SetAlgoLined
  $k=0$\\
  initial guess $\theta^0$ using method presented in Section~\ref{sec:thetaInit}\\
  compute $(H_1^0,H_2^0)$ from $\theta^0$ using Alg.~\ref{alg:H1H2reg}\\
  compute $\epsilon^0=E(\theta^0,H_1^0,H_2^0)$\\
  \While{$\epsilon^k < \epsilon^{k-1}$}{
    $k=k+1$\\
    compute $\theta^k$ from $(H_1^{k-1},H_2^{k-1})$ using Alg.~\ref{alg:theta}\\
    compute $(H_1^{k},H_2^{k})$ from $\theta^k$ using Alg.~\ref{alg:H1H2reg}\\
    compute $\epsilon^k=E(\theta^k,H_1^k,H_2^k)$
  }
\caption{Solver for $(\theta,H_1,H_2)$}
\label{alg:generalH1H2theta}
\end{algorithm}

\begin{figure}[h!t] 
  \begin{center}
    \includegraphics[width=0.29\textwidth]{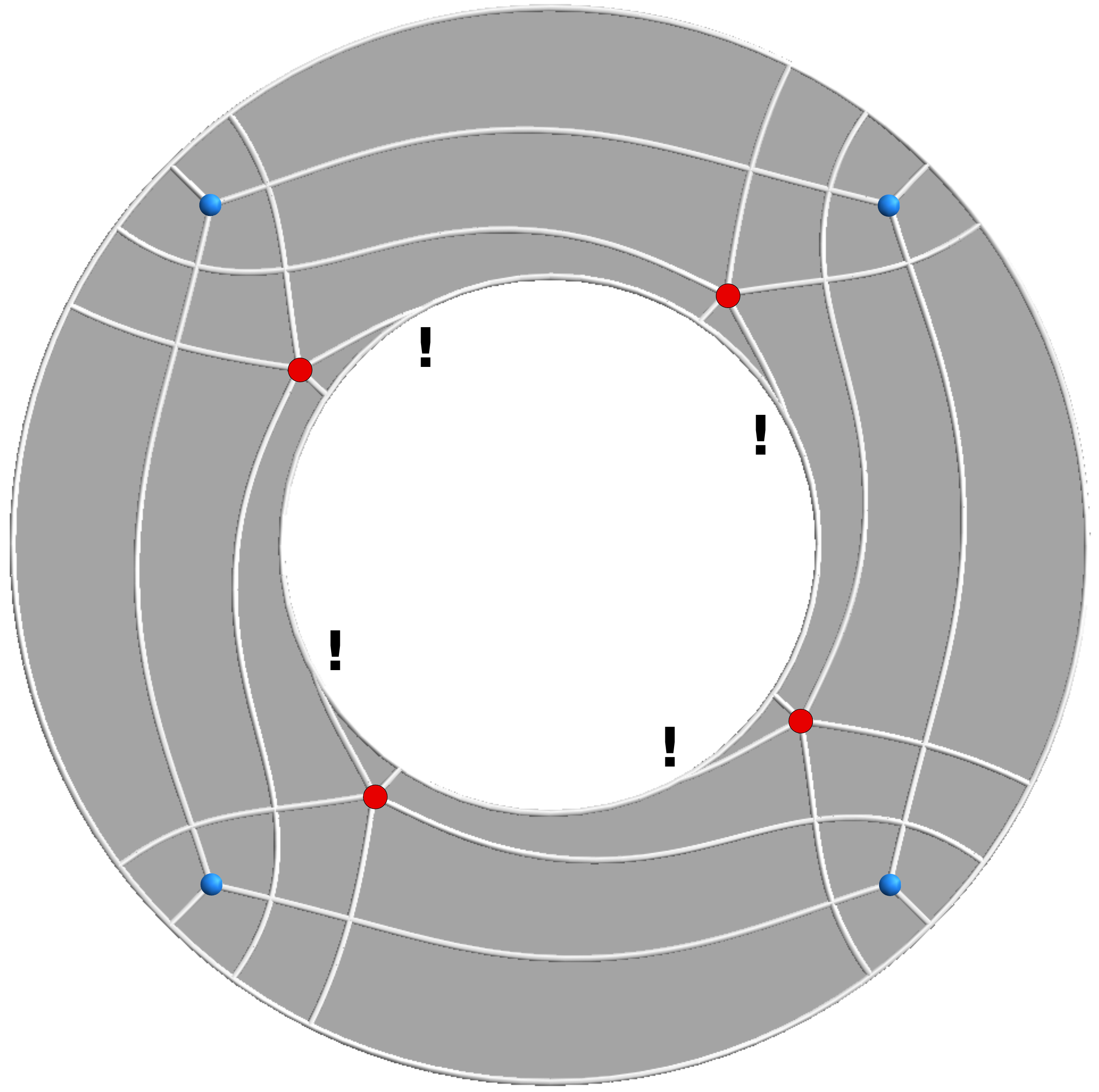}
    \includegraphics[width=0.29\textwidth]{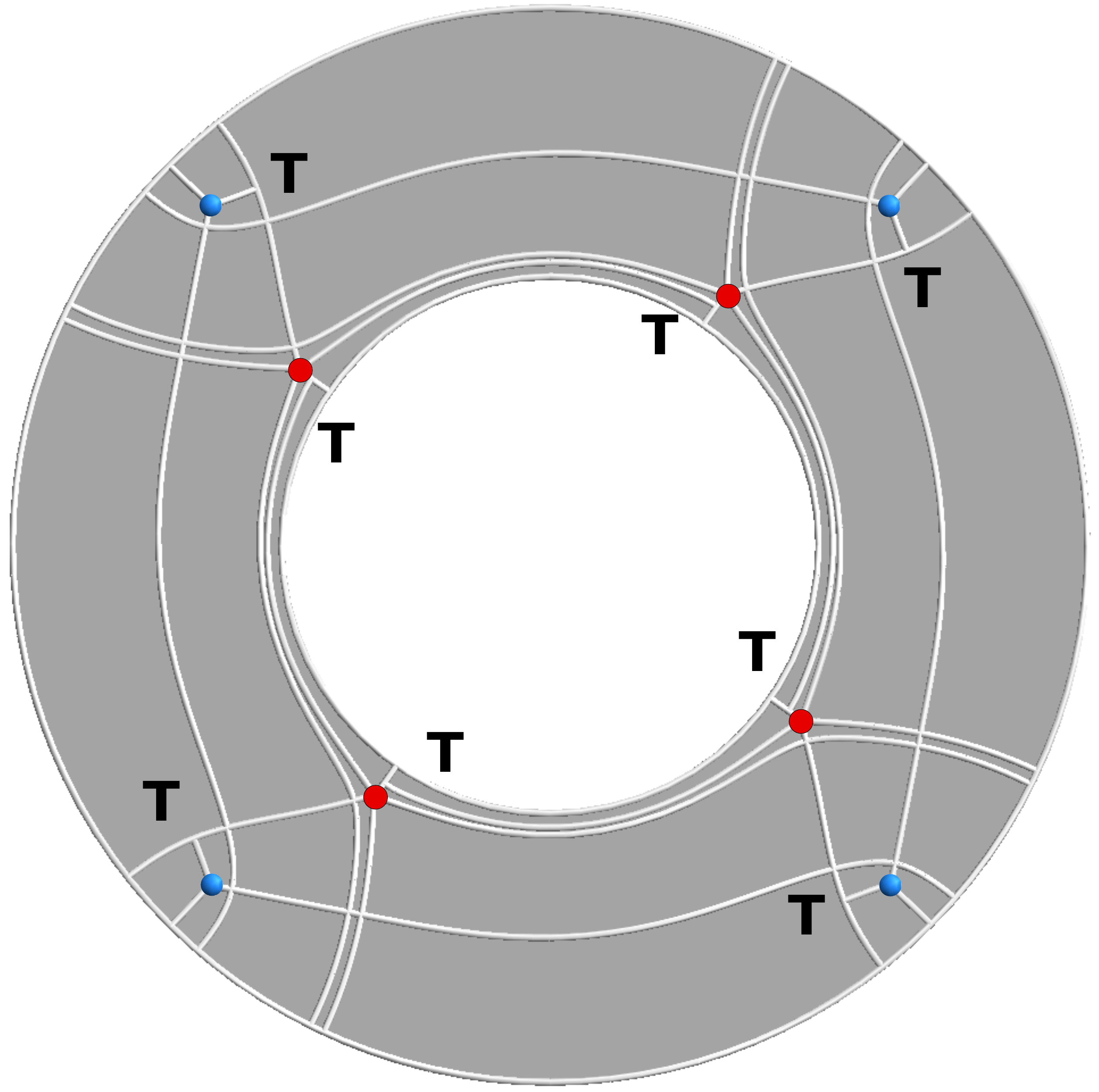}
    \caption{From left to right: quad layout obtained from a cross-field generated with dual $H$ function, and a quad layout obtained with the method presented in Section~\ref{sec:thetaInit}. The limit cycles are cut upon their first orthogonal intersection, therefore creating T-junctions.}
    \label{fig:ex1Int} 
  \end{center}
\end{figure}

For the sake of simplicity the motivational example, presented in Fig.~\ref{fig:ex1Int}, is planar and chosen to be topologically equivalent to a torus. A set of four of index $1$ and four of  index -$1$ singularities whose locations are not fulfilling the Abel-Jacobi condition is imposed. Consequently, a cross-field generated using the $H$ function will not be boundary aligned, and a cross-field generated using the method presented in Section~\ref{sec:thetaInit} will not be integrable and therefore will generate limit cycles.

\begin{figure}[h!t]
  \begin{center}
    \includegraphics[width=0.65\textwidth]{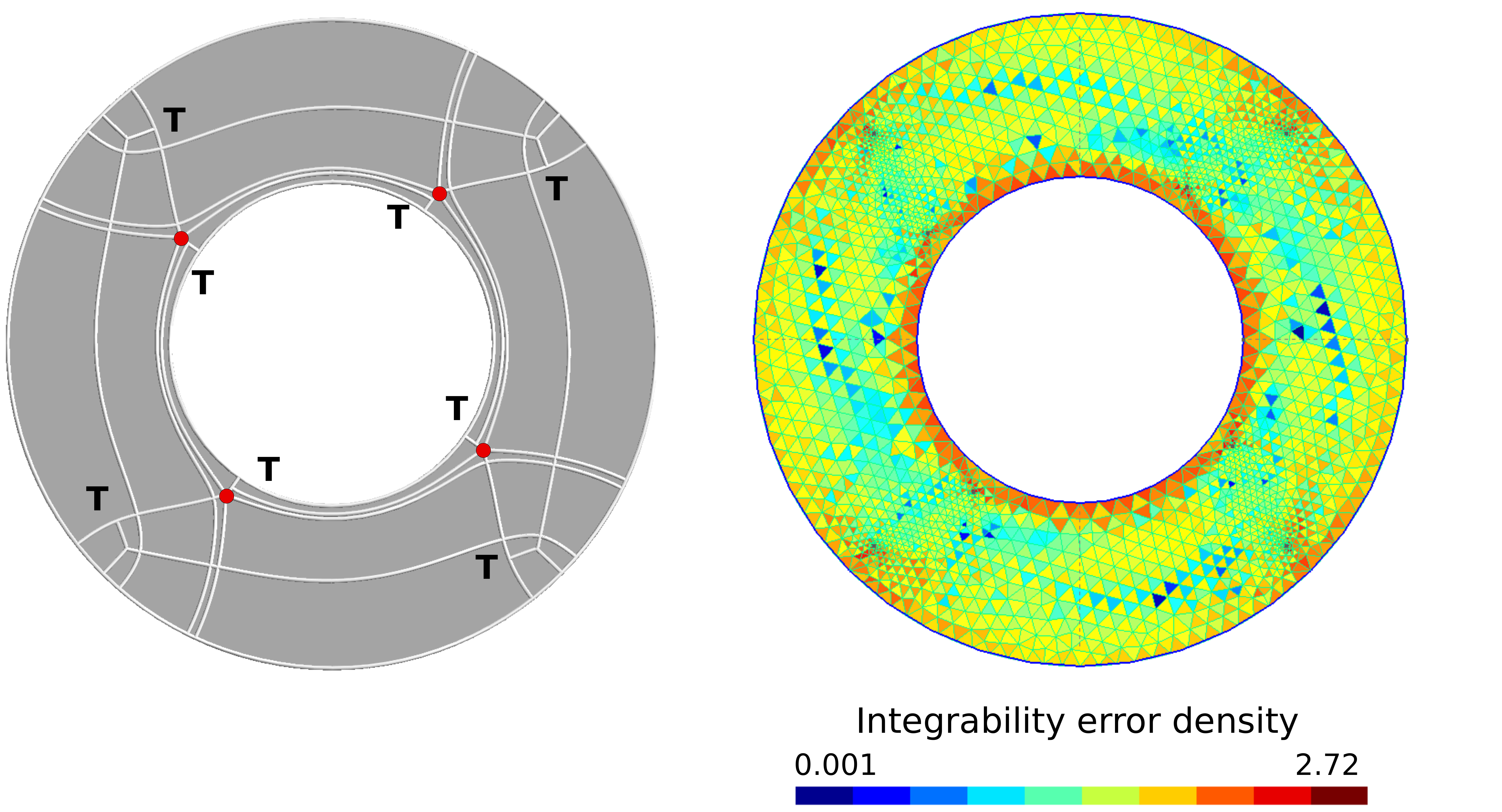}
    \caption{From left to right: quad layout obtained at initialization, and the integration error density on $\Omega$. The total integration error is $E=0.307898$.}
    \label{fig:ex1IntInit} 
  \end{center}
\end{figure}

\begin{figure}[h!t] 
  \begin{center}
    \includegraphics[width=0.65\textwidth]{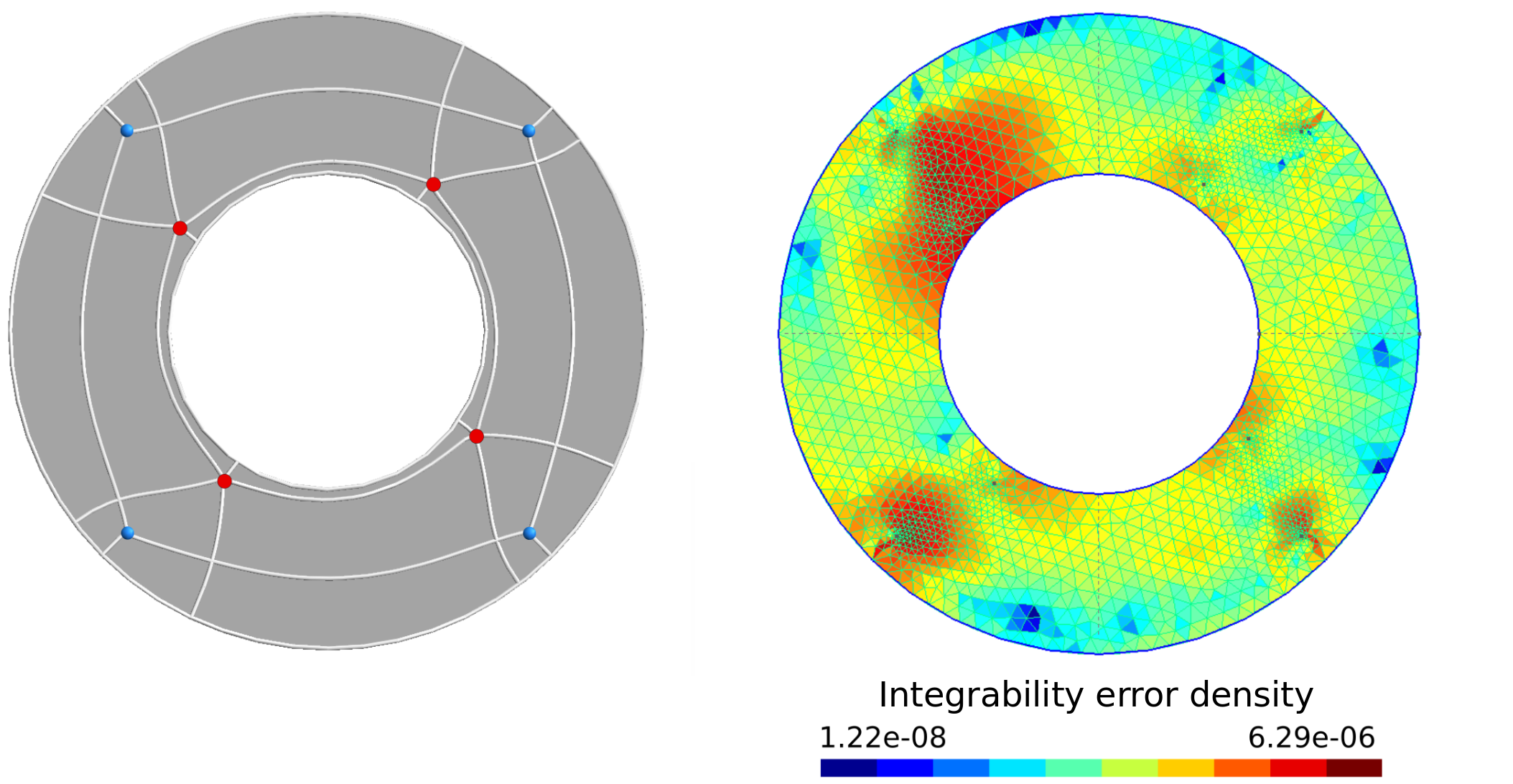}
    \caption{From left to right: quad layout obtained, and the integration error density on $\Omega$. The total integration error at convergence is $E=1.45639e-06$.}
    \label{fig:ex1IntConv} 
  \end{center}
\end{figure}

The method presented here is applied to compute an integrable boundary-aligned cross-field. Figure~\ref{fig:ex1IntInit} represents the cross-field used as an initial guess and Fig.~\ref{fig:ex1IntConv} is the one obtained at Algorithm~\ref{alg:generalH1H2theta} convergence.

Figure~\ref{fig:ex1IntInit} demonstrates that integrability error density is not 
concentrated in certain regions, but rather quite uniformly spread over the domain. This suggests that addressing the integrability issue cannot be performed via local modifications but only via the global one, \textit{i.e.}, the convergence of the presented non-linear problem. Figure~\ref{fig:ex1IntConv} shows that generating a limit cycle-free 2D cross-field can indeed be done by solving Eq.~(\ref{eq:chap3:intManifold2}). Nevertheless, this problem is highly non-linear and ill-defined, and solving it turns out to be difficult.

\begin{figure}[h!] 
  \begin{center}
    \includegraphics[width=0.29\textwidth]{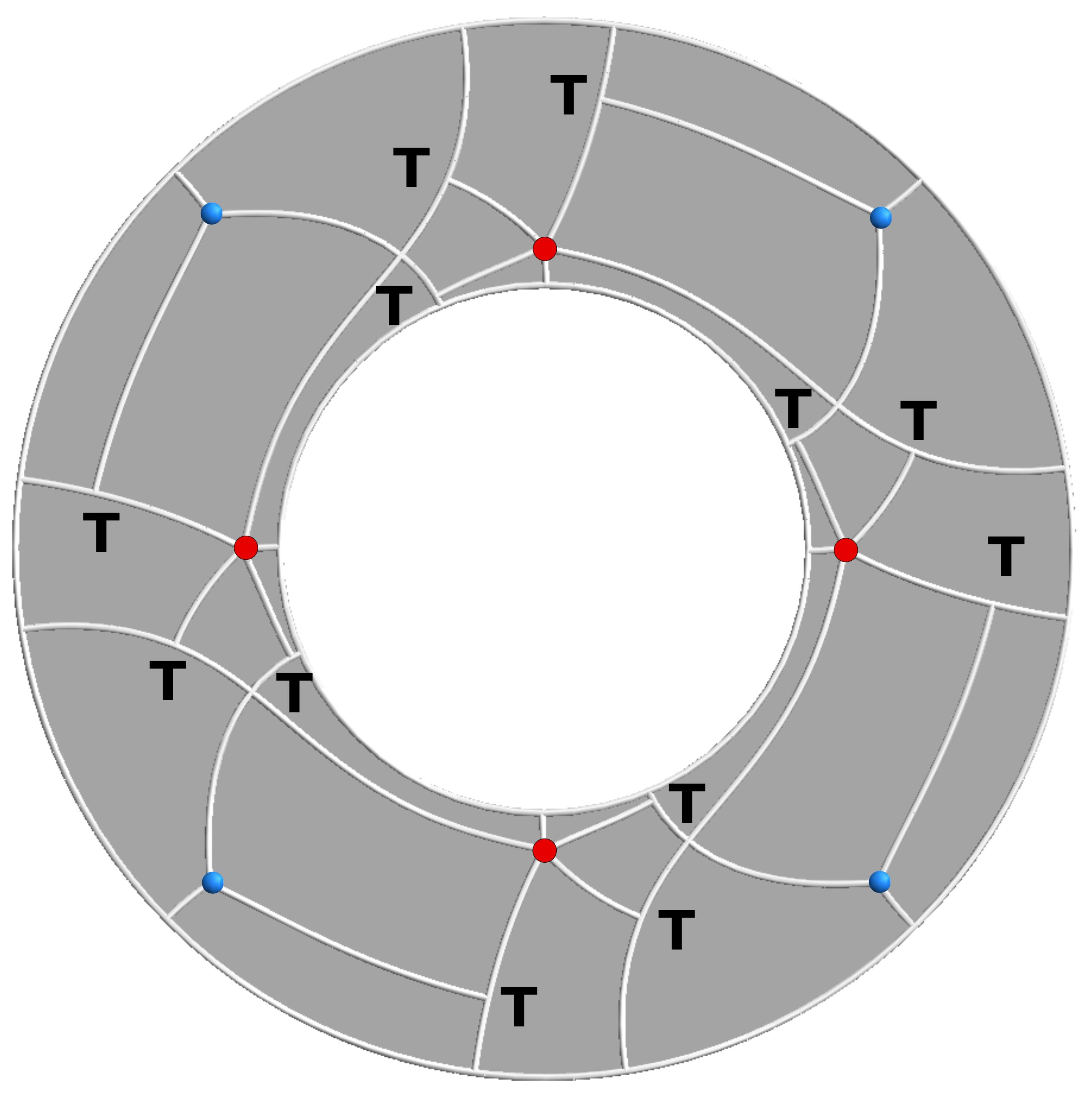}
    \includegraphics[width=0.29\textwidth]{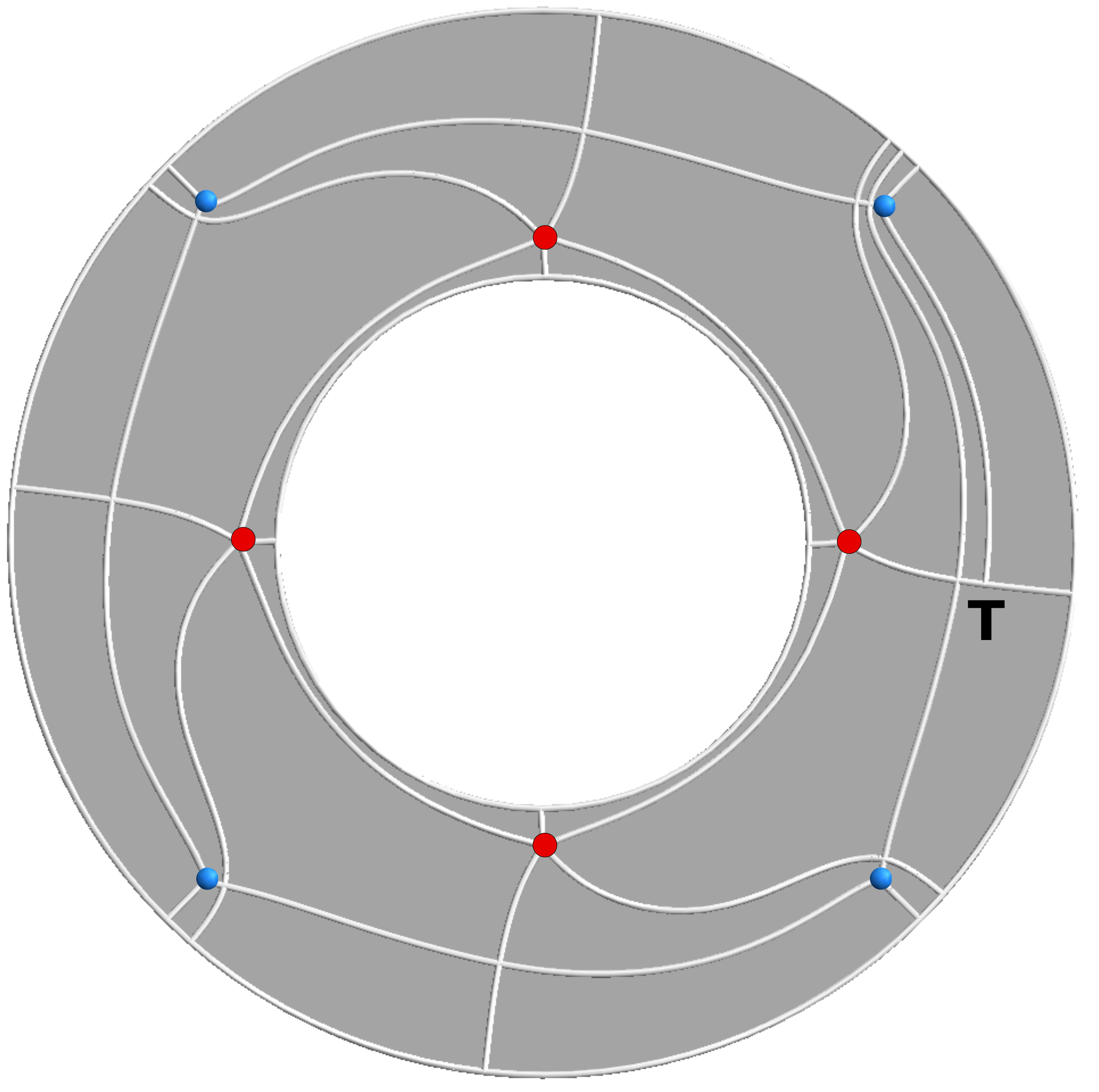}
    \caption{Left: quad layout obtained at initialization, the total integration error is $E=0.842169$. Right: quad layout obtained at convergence, the total integration error is $E=0.013597$}
    \label{fig:ex2IntConv} 
  \end{center}
\end{figure}

The method proposed here works well when initialization is not far from an integrable solution, \textit{i.e.}, when the imposed singularity set obeys Abel-Jacobi’s conditions. Otherwise, it does not converge up to the desired solution by reaching a local minimum $(\bar \theta,\bar H_1,\bar H_2)$ which does not satisfy $E(\bar \theta,\bar H_1,\bar H_2)=0$, as illustrated in Fig.~\ref{fig:ex2IntConv}. Although, it is interesting to note that, even without the presented method's convergence, the number of T-junctions dramatically decreases and the valid solution, in the opinion of authors, can be ``intuitively presumed''.

\vspace*{0.5cm}
\section{Conclusion and Future Work}
\label{S:Conclusion}

We presented the mathematical foundations for the generation of integrable cross-field on 2D manifolds based on user-imposed singularity configuration with both isotropic and anisotropic scaling.
Here, the mathematical setting is constrained by the Abel-Jacobi conditions for a valid singularity pattern. With the automatic algorithms to check and optimize singularity configuration (as recently presented in \cite{Chen:2019, Lei:2020, zheng2021quadrilateral}), the developed framework can be used to effectively generate both an isotropic and an anisotropic block-structured quad mesh with preserved singularity distribution. When it comes to computational costs of our cross-field generation, the formulation with isotropic scaling $H$ takes solving only two linear systems, and the anisotropic one $(H_1, H_2)$ represents a non-linear problem.

An attractive direction for future work includes, although it is not limited to, working with the user-imposed size map. By using the integrable cross-field formulation relying on two sizing fields $H_1$ and $H_2$, it would be possible to take into account the anisotropic size field to guide the cross-field generation. The size field obtained from the generated cross-field would not precisely match the one prescribed by the user, but it would be as close as possible to the singularity configuration chosen for the cross-field generation.

It is important to note that employing the presented framework in the $3$D volumetric domain would be possible only for a limited number of cases, in which the geometric and topological characteristics of the volume (more details in \cite{fogg2018singularities,white2000automatic}) allow the use of cross-field guided surface quad mesh for generating a hex mesh.

\newpage

\bibliographystyle{model1-num-names}
\bibliography{sample_arxiv.bib}
\end{document}